\documentclass{interact}
\usepackage[english]{babel}
\usepackage{amssymb, amsmath, amsthm} 
\usepackage{color}
\usepackage[font=small,labelfont=bf]{caption}
\usepackage{hyperref}
\usepackage{subfigure}

\numberwithin{equation}{section}

\usepackage{amsthm}
\makeatletter
\def\th@plain{%
  \thm@notefont{}
  \itshape 
}
\def\th@definition{%
  \thm@notefont{}
  \normalfont 
}
\makeatother

\newtheorem{thm}{Theorem}[section]
\newtheorem{prop}[thm]{Proposition}

\theoremstyle{definition}
\newtheorem{defn}[thm]{Definition}

\newtheorem{oss}[thm]{Remark}

\newcommand{\Z}{\mathbb{Z}}
\newcommand{\N}{\mathbb{N}}
\newcommand{\R}{\mathbb{R}}

\renewcommand{\H}{\mathcal{H}}
\newcommand{\I}{\mathcal{I}}
\newcommand{\D}{\mathcal{D}}

\renewcommand{\epsilon}{\varepsilon}

\def\e{\epsilon}

\DeclareMathOperator{\dist}{dist}

\begin{document}

\title{Motion of discrete interfaces on the triangular lattice}

\author{{\scshape Giovanni Scilla}\\
Department of Mathematics and Applications ``R. Caccioppoli''\\ University of Naples ``Federico II''\\
Via Cintia, Monte S. Angelo - 80126 Naples \\
(ITALY)}
\date{}

\maketitle

\begin{abstract}\noindent
We study the motion of discrete interfaces driven by ferromagnetic interactions on the two-dimensional triangular lattice by coupling the Almgren, Taylor and Wang minimizing movements approach and a discrete-to-continuum analysis, as introduced by Braides, Gelli and Novaga in the pioneering case of the square lattice. We examine the motion of origin-symmetric convex ``Wulff-like'' hexagons, i.e. origin-symmetric convex hexagons with sides having the same orientations as those of the 
hexagonal Wulff shape related to the density of the anisotropic perimeter 
$\Gamma$-limit of the ferromagnetic energies as the lattice spacing vanishes.
We compare the resulting limit motion 
with the corresponding evolution by crystalline curvature with natural mobility.
\end{abstract}


\begin{amscode}
35B27(Primary), 74Q10, 53C44, 49M25(Secondary).
\end{amscode}

\section{Introduction}
This paper is concerned with the variational motion of discrete interfaces arising from nearest-neighbours ferromagnetic-type interactions on the 2D triangular lattice. Our analysis aims to do a first step in the challenging and still largely open problem of characterizing the evolution of discrete interfacial energies driving more general atomistic systems in presence of dissipation. Indeed, the triangular lattice is the natural framework related to some discrete problems in crystallization (see, e.g.,~\cite{AYFS} and the references therein), fracture mechanics \cite{BG,BS2} and some physical models for two-dimensional fluids as the Bell-Lavis model~\cite{Bell1,Bell2}. Since crystalline perimeter energies can be approximated by lattice energies via $\Gamma$-convergence (see, e.g., \cite{ABC}) and arise in the study of evolutions by anisotropic curvature \cite{ATW83, CC, TCH1, TCH2, Ta0, Ta, Ta2}, the problem we address is also motivated by the analysis of the discreteness effects on such motions and their numerical approximation \cite{Chamb, ChambNov}. Moreover, our motions with underlying lattice can be interpreted as a simple version of the geometric evolutions in heterogeneous environments (see, e.g., \cite{Bat, Dirr}).\\ 
\noindent
{\bf The Almgren, Taylor and Wang approach to curvature driven motions.} The analysis will be carried over by using the minimizing-movements scheme of Almgren, Taylor and Wang~\cite{ATW83} for geometric evolutions driven by crystalline curvature. This consists in introducing a time scale $\tau$, an initial set $E_0^\tau$ and iteratively defining
a sequence of sets $\{E^\tau_k\}_{k\geq1}$ as minimizers of  
\begin{equation}
\min \left\{\int_{\partial E}\|\nu\|_1\mathrm{d}\H^1+{1\over \tau} \int_{E\triangle E^\tau_{k-1}}\dist(x,\partial E^\tau_{k-1})\,\mathrm{d}x\right\} ,
\label{ATW}
\end{equation}
among the sets of finite perimeter, where the first term is the \emph{crystalline perimeter} of $E$ 
and $\int_{E\triangle F}\dist(x,\partial F)\,\mathrm{d}x$ accounts for the $L^2$-distance between the boundaries of sets $E,F$; subsequently, computing a time-continuous limit (the \emph{flat flow}) $E(t)$ of the piecewise constant interpolations $E^\tau(t):=E^\tau_{\lfloor t/\tau\rfloor}$ 
as $\tau\to0$, which defines 
the motion by \emph{crystalline curvature}~\cite{AT95}
\begin{equation}
V(t)=\kappa(t),
\label{motcurv}
\end{equation}
as introduced independently by J. Taylor~\cite{Ta} and Angenent-Gurtin~\cite{AG}. According to \eqref{motcurv}, each side moves inward along its normal direction with a velocity $V$ coinciding with its crystalline curvature $\kappa$, that is proportional to the inverse of its length. In general, when dealing with 
an \emph{anisotropic perimeter} of the form
\begin{equation*}
\Phi(E)=\int_{\partial E}\varphi(\nu)\,\mathrm{d}\mathcal{H}^1,
\end{equation*}
where $\varphi$ is a norm on $\mathbb{R}^2$, {existence and uniqueness for the motion by crystalline curvature are simply proved for the class of ``good'' polygonal curves (see \cite[Prop. 2.1.1]{Ta}). That is the case, for instance,} if the initial set is a convex \emph{Wulff-like} set; i.e., it has a polygonal boundary 
whose sides have the same exterior unit normal vectors and form the same angles as those of the \emph{Wulff shape} $\mathcal{W}_\varphi$ of the density $\varphi$, where
\begin{equation}
\mathcal{W}_\varphi:=\Bigl\{x\in\mathbb{R}^2\,|\,\langle x,\nu\rangle\leq\varphi(\nu)\quad\text{ for every }\,\nu\mbox{ such that }|\nu|=1\Bigr\},
\end{equation}
$\langle\cdot,\cdot\rangle$ being the scalar product on $\R^2$. It is well known that $\mathcal{W}_\varphi$ is a centrally symmetric convex polygon and coincides with the unit ball $\{\varphi^\circ\leq1\}$ of the dual norm $\varphi^\circ$ (see, e.g., Morgan~\cite{FM}).

The Almgren, Taylor and Wang scheme implemented for energies
\begin{equation}
\int_{\partial E}\varphi(\nu)\,\mathrm{d}\mathcal{H}^1+\frac{1}{\tau}\int_{E\triangle E_{k-1}^\tau}\inf_{y\in\partial E_{k-1}^\tau}\varphi^\circ(x-y)\,\mathrm{d}x
\label{ATWgen}
\end{equation}
gives, in the limit as $\tau\to0$, the motion by crystalline curvature with \emph{natural mobility}
\begin{equation}
V(t)=M(n)\kappa(t),
\label{mobequa}
\end{equation}
where the mobility $M=\varphi$ is a function of the unit normal vector $n$ of the side (see, e.g., \cite{Ta0, TCH1, TCH2, Ta2, CC}). More precisely, the evolution is governed by equations
\begin{figure}[htbp]
\centering
\def\svgwidth{150pt}
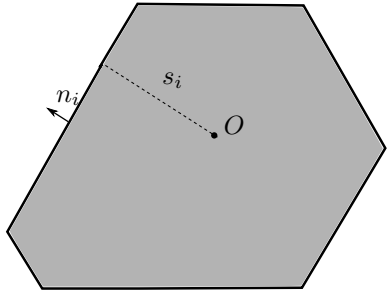
\caption{The function $s_i(t)$ is the distance from a fixed origin $O$ of a side with length $L_i(t)$ and normal vector $n_i$.}\label{fig:dist}
\end{figure}
\begin{equation}
\dot{s}_i(t)=-\varphi(n_i)\frac{\Lambda(n_i)}{L_i(t)},\quad i=1,\dots,m,
\label{equamotions}
\end{equation}
where $m$ is the number of sides of the initial set, $s_i$ is the distance from a fixed origin $O$ of a side with normal vector $n_i$ and length $L_i$, and $\Lambda(n_i)$ is the length of the side of the Wulff-shape having $n_i$ as normal vector (see Fig.~\ref{fig:dist}).

For example, if $\varphi(\nu)=\|\nu\|_1$ then its Wulff shape is the coordinate square $\mathcal{W}_{\|\cdot\|_1}=\{\|x\|_{\infty}\leq1\}$, 
$n_i\in\{\pm e_1, \pm e_2\}$, $i=1,2,3,4$, $\varphi(n_i)\equiv1$ and $\Lambda(n_i)\equiv2$.
Thus, in the case of initial datum a coordinate rectangle (a convex `Wulff-like' set), the evolution by crystalline curvature \eqref{equamotions} is a rectangle with the same centre and sides of lengths $L_1(t), L_2(t)$ solution to the system of ordinary differential equations
\begin{equation}
\begin{cases}\displaystyle \dot L_1(t)= -{4\over L_2(t)}\cr\cr
\displaystyle \dot L_2(t)= -{4\over L_1(t)}, \quad t\in[0,T].\end{cases}
\label{crystarect}
\end{equation}
\noindent
{\bf Motion of discrete interfaces: the state of the art.} In a forerunner paper by Braides, Gelli and Novaga~\cite{BGN} the Almgren, Taylor and Wang approach (\ref{ATW}) has been coupled with a homogenization procedure via $\Gamma$-convergence. In this case the perimeters and the distances depend on a small parameter $\e>0$ (the space scale), and consequently, after introducing a time scale $\tau$, $E_0^{\tau,\epsilon}$ such that ${\rm d}_\mathcal{H}(E_0^{\tau,\epsilon},E_0)\to0$ for some regular limit set $E_0$, the time-discrete motions are the sets $E^{\tau,\e}_k$ defined iteratively by 
\begin{equation}
E^{\tau,\e}_k \in \text{argmin} \Bigl\{ P_\e(E)+{1\over \tau} \int_{E\triangle E_{k-1}^{\tau,\epsilon}} d_\epsilon(x,\partial E^{\tau,\e}_{k-1})\,\mathrm{d}x\Bigr\},\quad k\geq1,
\label{discrtrajec}
\end{equation}
where the minimization is over finite unions of squares with side-length $\varepsilon$.
The energies $P_\e$ are {\em discrete ferromagnetic energies} (see Remark~\ref{justify} for a physical interpretation), defined on subsets of the square lattice $E\subset \e\Z^2$, of the form
\begin{equation}
P_\e(E)= \frac{1}{2}\,\e\, \#\Bigl\{(i,j)\in\e\Z^2\times\e\Z^2: i\in E, j\not\in E, \ |i-j|=\e\Bigr\},
\label{ene}
\end{equation}
\\
each couple $(i,j)$ being accounted twice. After a piecewise constant identification of $E$ with a subset of $\mathbb{R}^2$, the $P_\varepsilon$ can be interpreted as the perimeter of $E$. The continuum ($\Gamma$\hbox{-})limit of these energies as $\epsilon\to0$ is the crystalline perimeter 
\begin{equation}
P(E)=\int_{\partial E}\|\nu\|_1\mathrm{d}\H^1,
\label{crysperimeter}
\end{equation}
as proved by Alicandro, Braides and Cicalese~\cite{ABC}. The distance $d_\epsilon$ in \eqref{discrtrajec} is a suitable discretization of the $L^\infty$-distance from the boundary of the previous set.

The time-continuous limit $E(t)$ of $E^{\tau,\e}_{\lfloor t/\tau\rfloor}$ defined in (\ref{discrtrajec}) then may depend how mutually $\e$ and $\tau$ tend to $0$ (Braides~\cite[Theorem~8.1]{Bra13}). In particular, if $\tau/\e\to0$ the limit motion will be \emph{pinned}; i.e., $E(t)\equiv E_0$ (in a sense, we can pass to the limit in $\tau$ first, and then apply the Almgren-Taylor-Wang approach). On the contrary, if $\tau/\e\to+\infty$ then the limit $E(t)$ will be 
the crystalline evolution related to the limit $P$ defined in \eqref{crysperimeter} (again, in a sense, in this case we can pass to the limit in $\e$ first). 
The relevant regime which gives the most information about all the limit evolutions is when $\tau/\e\to\gamma\in(0,+\infty)$. In this case, when the initial datum is a coordinate rectangle, the resulting evolution is still a coordinate rectangle and, in case of uniqueness, the side-lengths $L_1(t), L_2(t)$ of this rectangle solve a system of ``degenerate'' ordinary differential equations
\begin{equation}
\begin{cases}\displaystyle \dot L_1(t)= -{2\over \gamma}\left\lfloor{2\gamma\over L_2(t)}\right\rfloor\\ \cr
\displaystyle \dot L_2(t)= -{2\over \gamma}\left\lfloor{2\gamma\over L_1(t)}\right\rfloor,\end{cases}
\label{sistema}
\end{equation}
for almost every $t$ until the extinction time. We note that the discontinuous form of the right-hand sides highlights that the microscopic motion is obtained by overcoming some energy barriers in a `quantized' manner. In particular, we have `pinning' of large rectangles: if both initial side-lengths are above the {\em pinning threshold} $\widetilde L=2\gamma$ then the right-hand sides in (\ref{sistema}) are zero and the motion is pinned. The limit cases of total pinning and continuous crystalline flow \eqref{crystarect} correspond to the limit values $\gamma=0$ and $\gamma=+\infty$, respectively. 

These unexpected features of the limit motions led many authors to investigate the sensibility of such evolutions to microstructure. Indeed, they are very sensitive and may depend on microscopic properties not detected in the limit description, as showed, e.g., by Braides and Scilla~\cite{BraSci} in case of periodic media and by Scilla~\cite{Sci14} in case of `low-contrast' periodic media; therein the dependence of the limit velocity on the curvature is described by a homogenized formula quite different with respect to \cite{BGN}. A random counterpart of the low-contrast setting has been provided by Ruf~\cite{Ruf}.  Recently, Braides, Cicalese and Yip~\cite{BCY} investigated the case of antiferromagnetic energies, in particular anti-phase boundaries between striped patterns, showing the appearance of some non-local curvature dependence velocity law reflecting the creation of some defect structure on the interface at the microscopic level. Braides and Solci~\cite{BraSo15}, instead, treated the motion through \emph{mushy layers} in high-contrast spin systems, dealing with the creation of bulk microstructure. An approach to time-reversed motions can be found in Braides and Scilla~\cite{BraSci2}, where a suitably scaled discrete version of the Almgren, Taylor and Wang scheme has been used with a negative perimeter term, thus forcing the minimizers to have a ``checkerboard'' structure. We mention also \cite{BMN, MN} for effective crystalline evolutions resulting from microstructures modeled through periodic forcing terms.

\noindent
{\bf Our result.} In this paper we perform a discrete-to-continuum analysis as in \cite{BGN} on a different lattice by choosing a suitable dissipation term in the energy, and we compare the resulting limit evolution with the corresponding crystalline motion related to the limit perimeter of the lattice energies. The main result is that, for convex initial sets whose geometry is compatible with the underlying lattice (i.e., ``Wulff-like'' hexagons), the limit evolution is still hexagonal but the shrinking velocity is slower.

For this, we consider the same lattice energies as in (\ref{ene}) labeled by the nodes of the \emph{triangular lattice} $\mathbb{T}$, generated, e.g., by vectors $\eta_1=(1,0)$ and $\eta_2=(1/2,\sqrt{3}/2)$.
A computation obtained by adapting the arguments in~\cite{ABC,BG} shows that the continuous limit as $\epsilon\to0$ of the energies (\ref{ene}) is the {anisotropic perimeter}
\begin{equation}
{P}(E)=\int_{\partial E}\varphi_{hex}(\nu)\mathrm{d}\H^1,
\label{hexagonalper}
\end{equation}
whose density $\varphi_{hex}$ is a norm 
with hexagonal symmetries, and the corresponding Wulff shape $\mathcal{W}_{hex}$ is a regular hexagon.

In the discrete formulation of \eqref{ATWgen}, we restrict ourselves to initial limit sets which are convex origin-symmetric {Wulff-like} hexagons (see Definition~\ref{wshlikeset} for a precise definition). These sets can be seen as the equivalent of coordinate rectangles in \cite{BGN}, and the motion of more general sets can be reduced to the study of these ones.

In the case of initial datum a convex Wulff-like hexagon, with Proposition~\ref{rectangleprop} we prove that the resulting evolution is a set of the same type. Indeed, each connected component of the evolution is a convex Wulff-like hexagon, since a Wulff-like convexification provides a competitor with less energy in the minimization problem. Then, an argument based on suitable translations towards any incenter of the previous set shows that the evolution is actually connected. If, in addition, the initial set is origin-symmetric, then the evolution preserves this property as well. Furthermore, the underlying lattice forces the velocities to be quantized. Indeed, each side of length $L_{i}(t)$ with exterior unit normal vector $n_i$ moves inward 
in the direction $n_i$, and its distance $s_i(t)$ from the origin reduces with velocity $v_{i}(t)$ satisfying the inclusions

\begin{equation*}
v_{i}(t)
\begin{cases}
=\displaystyle\frac{\sqrt3}{2\gamma}\left\lfloor\frac{\alpha_{hex}\gamma}{L_i(t)}\right\rfloor, & \mbox{if \, $\displaystyle\frac{\alpha_{hex}\gamma}{L_i(t)}\not\in\mathbb{N}$}\\
\\
\in\displaystyle\frac{\sqrt3}{2\gamma}\left[\displaystyle\left(\frac{\alpha_{hex}\gamma}{L_i(t)}-1\right),\displaystyle\frac{\alpha_{hex}\gamma}{L_i(t)}\right], & \mbox{if \, $\displaystyle\frac{\alpha_{hex}\gamma}{L_i(t)}\in\mathbb{N}$},
\end{cases}
\end{equation*}
where the mobility factor $\alpha_{hex}$ is equal to $\frac{16}{9}$.


As a consequence, in the case of a unique evolution, the distances $s_i(t), i=1,\dots,6$ 
solve the system of degenerate ordinary differential equations
\begin{equation}\label{systems}
\dot{s}_i(t)=-{\sqrt3\over 2\gamma}\left\lfloor\frac{\alpha_{hex}\gamma}{L_{i}(t)}\right\rfloor, \quad i=1,\dots,6
\end{equation}
for almost every $t$. For each side of length $L_i(t)$, \eqref{systems} translates into
\begin{equation}\label{system3}
\dot{L}_i(t)=\frac{1}{\gamma}\left\lfloor\frac{\alpha_{hex}\gamma}{L_{i}(t)}\right\rfloor-\frac{1}{\gamma}\left(\left\lfloor\frac{\alpha_{hex}\gamma}{L_{i-1}(t)}\right\rfloor+\left\lfloor\frac{\alpha_{hex}\gamma}{L_{i+1}(t)}\right\rfloor\right),\quad i=1,\dots,6,
\end{equation}
where the labelling of the sides is intended to be modulo 6.
%
Clearly, by \eqref{system3}, we have pinning of large hexagons: if the side lengths $L_i^0$ of the initial limit set comply with the condition $\displaystyle\min_{1\leq i\leq6}\{L_i^0\}>\displaystyle\alpha_{hex}\gamma$, then none of the sides can move and the motion coincides identically with the initial set. The pinning threshold is computed by imposing that the minimal inward displacement of a side along its normal direction is not energetically convenient (see Section~\ref{pinning}). We point out that also some \emph{partial pinning} phenomena may happen; that is,  
a side stays pinned until it shortens sufficiently due to the motion of the adjacent sides (see Remark~\ref{esempiopartial} for an example).

We note that in the particular case when the initial datum is a (sufficiently small) regular Wulff-like hexagon (\emph{motion of the Wulff shape}), $s_i(t)\equiv s(t)$ is independent of the sides and coincides with the apothem. In this case we have a self-similar evolution and the system (\ref{system3}) reduces to the single equation
\begin{equation}\label{system2}
\dot{L}(t)=-{1\over \gamma}\left\lfloor\frac{\alpha_{hex}\gamma}{L(t)}\right\rfloor,
\end{equation}
while the corresponding crystalline evolution with natural mobility, according to \eqref{equamotions}, is given by
\begin{equation}
\dot{L}(t)=-\alpha_{hex}\frac{1}{L(t)}.
\label{corrcrys}
\end{equation}
This shows that the main effect of discreteness on the limit motion is to slow down the corresponding crystalline evolution (\ref{corrcrys}), that, however, can be retrieved from (\ref{system2}) in the limit as $\gamma\to+\infty$.

\noindent
{\bf Plan of the paper.} The paper is organized as follows. In Section~\ref{setting} we fix notation and introduce the energies we will consider on the triangular lattice. We then formulate the discrete analogous of the Almgren, Taylor and Wang scheme \eqref{ATWgen}, according to~\cite{BGN}. Section~\ref{hexagonal} contains the proof of the convergence of the discrete scheme in the case of an origin-symmetric convex Wulff-like initial set. 
In Section~\ref{pinning} we compute the pinning threshold and the description of the  limit motion is contained in Section~\ref{limitmotion}, where we also compare it with the corresponding crystalline evolution. 

\section{Setting of the problem}\label{setting}

If $x=(x^1,x^2)\in \R^2$ we set $\|x\|_1=|x^1|+|x^2|$, $\|x\|_\infty=\max\{|x^1|,|x^2|\}$ and the usual scalar product in $\R^2$ will be denoted by $\langle\cdot,\cdot\rangle$. If $A$ is a Le\-bes\-gue\hbox{-}measurable set we denote by $|A|$ its two\hbox{-}dimensional Lebesgue measure. The symmetric difference between two sets $A$ and $B$ in $\mathbb{R}^2$ is denoted by $A\triangle B$, their \emph{Hausdorff distance} is defined by
\begin{equation*}
{\rm d}_\mathcal{H}(A,B)=\max\left\{\sup_{a\in A}\dist(a,B), \sup_{b\in B}\dist(b,A)\right\},
\end{equation*}
where $\dist(x,E)$ is the distance of the point $x$ to the set $E$ defined, as usual, by $\dist(x,E)=\inf_{y\in E}|x-y|$. We say that a sequence of sets $\{E_\epsilon\}$ converges to $E$ in the Hausdorff sense as $\epsilon\to0$ if and only if ${\rm d}_{\mathcal H}(E_\epsilon,E)\to0$ and $E$ is closed.

If $E$ is a set of finite perimeter then $\partial^*E$ is its reduced boundary (see, for example~\cite{Bra98}) and the 1-dimensional Hausdorff measure of $\partial^*A$ is denoted by $\mathcal{H}^1(\partial^*A)$. The 
measure-theoretical inner normal to $E$ at a point $x$ in $\partial^*E$ is denoted by $\nu=\nu_E(x).$

If $e=(e^1,e^2)$ is a vector, then we denote by $e^\perp$ the anticlockwise rotation of $\pi/2$ of $e$, that is, $e^\perp=(-e^2, e^1)$. If $\{a_i\}, i=1,\dots,N$ is a finite set of vectors, then $conv(a_1,a_2,\dots,a_N)$ is the \emph{convex hull} of vectors $\{a_i\}$.

Let $\varphi$ be a norm on $\R^2$ and let $A\subset\R^2$ be a convex set. The \emph{inradius} $r_A$ of $A$ is the radius of the largest balls contained in $A$ and can be defined as
\begin{equation*}
r_A:=\max_{x\in A}\min_{y\in\partial A}\varphi(x-y).
\end{equation*} 
Correspondingly, the centers of the balls of largest radius $r_A$ inscribed in $A$ are called the \emph{incenters of $A$}, and we denote by $\mathcal{M}(A)$ their set:

\begin{equation}
\mathcal{M}(A):=\Bigl\{x\in A:\, \min_{y\in\partial A}\varphi(x-y)=r_A\Bigr\}.
\label{kernel}
\end{equation}
We note that, in general, $\mathcal{M}(A)$ is not a singleton; it can be shown (see, e.g.,~\cite[Lemma~3.3]{BraMai}) that $\mathcal{M}(A)$ is a closed convex set with $|\mathcal{M}(A)|=0$, and its dimension is at most 1 (i.e., $\mathcal{M}(A)$ is either a point or a line segment).

\subsection{Ferromagnetic energies on triangular lattice}
\begin{figure}[htbp]
\centering
\def\svgwidth{120pt}
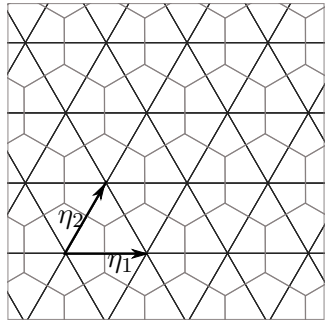
\caption{Triangular lattice and its dual (hexagonal) lattice.}\label{fig:0}
\end{figure}

We consider the \emph{triangular lattice} $\mathbb{T}=\{a\eta_1+b\eta_2|\, a,b\in\mathbb{Z}\}$, where $\eta_1=(1,0)$ and $\eta_2=(1/2,\sqrt{3}/2)$ (see Fig.~\ref{fig:0}). 
We also define $\eta_3:=\eta_1-\eta_2$, $\mathcal{S}:=\{\pm \eta_1,\pm\eta_2,\pm\eta_3\}$ the set of the \emph{unitary vectors} in the lattice $\mathbb{T}$ and $\mathcal{R}:=\{\eta^\perp:\, \eta\in\mathcal{S}\}$ the set of \emph{coordinate directions}. With fixed $\epsilon>0$, we introduce interfacial energies defined on subsets $\mathcal{I}\subset\epsilon\mathbb{T}$ as 
\begin{equation}
P_\e(\mathcal{I})= \,\frac{\sqrt 3}{3}\left(\frac{1}{2}\e\, \#\Bigl\{(i,j)\in\e\mathbb{T}\times\e\mathbb{T}: i\in \I, j\not\in \I, \ |i-j|=\e\Bigr\}\right),
\label{energies}
\end{equation}
where each couple $(i,j)$ is accounted twice. 

To study the continuous limit as $\e\to0$ of these energies it is customary to identify each subset of $\e\mathbb{T}$ with a measurable subset of $\R^2$, in such a way that equi-boundedness of the energies implies pre-compactness of such sets in the sense of sets of finite perimeter. This identification is as follows: we denote by $H$ the cell of the dual lattice of $\mathbb{T}$ (see Fig~\ref{fig:0}); that is, the closed regular hexagon of center $0$ and side-length $\sqrt 3/3$ defined by 
$H=\frac{\sqrt 3}{3}conv(\mathcal{R})$. For every ${i}\in\epsilon\mathbb{T}$, we denote by $H_\epsilon( {i})= {i}+\epsilon H$ the closed regular hexagon with side-length $\frac{\sqrt 3}{3}\epsilon$ and centered in $ {i}$. We refer to each side of $H_\epsilon( {i})$ as a \emph{cell side}. 

To a set of indices $\I\subset\epsilon\mathbb{T}$ we associate the set

\begin{equation*}
E_{\I}=\bigcup_{ {i}\in \I}H_\epsilon( {i})\subset\R^2.
\end{equation*} 
\\
The space of \emph{admissible sets} related to indices in the two\hbox{-}dimensional triangular lattice is then defined by
\begin{equation}
\D_\epsilon:=\Bigl\{E\subseteq\R^2:\quad E=E_\I\text{ for some $\I\subseteq\epsilon\mathbb{T}$}\Bigr\},
\label{admissibleset}
\end{equation}
and for each $E=E_\I\in \D_\epsilon$ we denote, with abuse of notation,
\begin{equation}\label{pabei}
P_\epsilon(E)=P_\epsilon(\I)=\mathcal{H}^1(\partial E).
\end{equation}

\begin{oss}\label{justify}
In order to justify the name of \emph{ferromagnetic energies}, we remark that the $P_\epsilon$ as in \eqref{energies} can be viewed as lattice energies; that is, depending on a discrete variable $u=\{u_i\}$ indexed by the nodes $i$ of $\e\mathbb{T}$, of the form
\begin{equation}
P_\epsilon(u)=\frac{\sqrt{3}}{3}\times\frac{1}{4}\sum_{|i-j|=\epsilon}\epsilon(1-u_iu_j),
\end{equation}
where $u_i=u(i)$ takes only the two values +1 and -1 (\emph{ferromagnetic energy for Ising spin system}). Here the factor $1/4$ is due to the fact that each couple $(i,j)$ of nearest neighbors is accounted twice. After
identifying $u$ with the set $E$ obtained as the union of all closed hexagons with side length $\frac{\sqrt{3}}{3}\epsilon$ and centers $i$
such that $u_i = 1$, the energy $P_\epsilon$ can be equivalently rewritten as a perimeter functional as in \eqref{pabei}, and hence can be interpreted as a \emph{discrete interfacial energy}.
\end{oss}

In~\cite{ABC} it is shown that the $\Gamma$-limit's domain of energies $P_\varepsilon$ is the family of sets of finite perimeter and its general form is
$$
\Phi(E)=\int_{\partial^*E}\varphi(\nu)\mathrm{d}\H^1,
$$
with $\varphi$ a convex, positively homogeneous of degree one function reflecting the symmetries of the underlying lattice.

%

We state without proof 
the $\Gamma\hbox{-}$convergence result for energies \eqref{energies}, noting that the argument may be deduced from a more general proof developed for vector-valued Lennard-Jones type interactions~\cite[Proposition 4.6]{BG}.

\begin{thm} [$\Gamma$\hbox{-}convergence of perimeter energies]\label{thmconv}
The energies $\text{P}_\epsilon$ defined by
\begin{equation}
\text{P}_\epsilon(E)=
\begin{cases}
\text{P}_\epsilon(\I), &\text{if }E=E_{\I}\in\mathcal{D}_\epsilon\\
+\infty, & \text{otherwise}
\end{cases}
\end{equation}
$\Gamma$\hbox{-}converge, as $\epsilon\to0$, with respect to the $L^1$-topology to the anisotropic perimeter functional
\begin{equation}
P(E)=\int_{\partial^* E}\varphi_{hex}(\nu)\,\mathrm{d}\H^1,
\label{functionalhexa}
\end{equation}
whose density $\varphi_{hex}$ is defined as 
\begin{equation}
\varphi_{hex}(\nu):=\frac{2}{{3}}\sum_{k=1}^3|\langle \nu, \eta_k\rangle|,
\label{minnorm}
\end{equation}
where $\eta_k$, $k=1,2,3$, are the unitary vectors of the lattice.
\end{thm}

\begin{oss}
If $\varphi_{hex}$ is defined as in \eqref{minnorm}, then
\begin{equation}
\varphi_{hex}(\nu)=\frac{4}{3}\max_{k=1,2,3}|\langle \nu, \eta_k\rangle|,\quad \mbox{ for every $\nu\in\R^2$,}
\end{equation}
and $\{\varphi_{hex}\leq1\}=\displaystyle\frac{\sqrt{3}}{2}conv(\pm \eta_1^\perp,\pm \eta_2^\perp, \pm \eta_3^\perp)$.
\end{oss}

\subsection{The dissipation term and the minimization scheme}\label{timemin} 

The choice of the dissipation term in the Almgren, Taylor and Wang scheme affects the mobility of the limit interface (see, e.g., \cite[Section~1]{Ta0} for a discussion). For instance, considering there the distance induced by the dual norm of the density of the perimeter term as in \eqref{ATWgen}, in the limit as $\tau\to0$ one retrieves the motion by crystalline curvature with natural mobility, governed by equations \eqref{equamotions}. Although this situation is not very general, as a fact of interest in this case the evolution of the Wulff shape is explicit and self-similar (see, e.g., \cite{BePa, CC}). Another motivation, purely practical, is that the level sets of the resulting distance have the symmetries of the Wulff shape, thus simplifying many computations.

Therefore, in order to define the dissipation term in the Almgren, Taylor and Wang scheme, we first notice that, since the linear function $\langle \xi,\nu\rangle$ on the nonempty compact convex polygon $\{\varphi_{hex}(\nu)\leq1\}$ attains its maximum at a vertex of the polygon (see, e.g., \cite[Corollary 32.3.2]{Rock}), the dual norm $\varphi^\circ_{hex}$ of $\varphi_{hex}$ is given by
\begin{equation}
\varphi_{hex}^\circ(\xi)=\frac{\sqrt{3}}{2}\max_{k=1,2,3}|\langle \xi, \eta_k^\perp\rangle|,\quad \xi=(\xi^1,\xi^2).
\label{duale}
\end{equation}
{Moreover, \eqref{duale} complies with $\varphi_{hex}^\circ(\mathbb{T}\backslash\{0\})=\frac{3}{4}\N$, since if $\xi\in\mathbb{T}\backslash\{0\}$, $\xi=n\eta_1+m\eta_2$ for some $n,m\in\Z$, then $\varphi_{hex}^\circ(\xi)=\frac{3}{4}\max\{|m-n|,|m+n|\}$.}\\

We define a notion of \emph{discrete distance} $d_{\varphi^\circ_{hex}}^\epsilon$ induced by $\varphi^\circ_{hex}$ as
\begin{enumerate}
\item[(1)] $d_{\varphi^\circ_{hex}}^\epsilon(H_\epsilon(i),H_\epsilon(j))=\varphi^\circ_{hex}(i-j)$,\,\,for every $i,j\in\varepsilon\mathbb{T}$;
\item[(2)] $d_{\varphi^\circ_{hex}}^\epsilon(E,F)=\inf\left\{d_{\varphi^\circ_{hex}}^\epsilon(H_\epsilon(i),H_\epsilon(j)):\,\,H_\epsilon(i)\in E,H_\epsilon(j)\in F \right\}$, for every $E,F\in\mathcal{D}_\epsilon$.
\end{enumerate}
Moreover, for every $x\in E$, we set
\begin{equation}
d_{\varphi^\circ_{hex}}^\epsilon(x,\partial F):=
\begin{cases}
\inf\left\{d_{\varphi^\circ_{hex}}^\epsilon(H_\epsilon(i),H_\epsilon(j)):\,\,x\in H_\epsilon(i), \, H_\epsilon(j)\in F \right\},\quad {\rm if}\,\, x\not\in F,\\
\inf\left\{d_{\varphi^\circ_{hex}}^\epsilon(H_\epsilon(i),H_\epsilon(j)):\,\,x\in H_\epsilon(i), \, H_\epsilon(j)\not\in F \right\},\quad {\rm if}\,\, x\in F.
\end{cases}
\label{discretedistance}
\end{equation}
Note that this is well defined as a measurable function, since its definition is unique outside the union of the boundaries of the hexagons $H_\epsilon$ (that are a negligible set).\\


Now, we introduce the same discrete minimization scheme as in \cite{BGN}. We fix a time step $\tau>0$ and define a discrete motion with underlying time step $\tau$ obtained by successive minimization. At each time step we will minimize the energy $\mathcal{F}_{\tau,\epsilon}:\D_\epsilon\times\D_\epsilon\to\R$ given by

\begin{equation}
\begin{split}
\mathcal{F}_{\tau,\epsilon}(E,F)&= P_\epsilon(E)+\frac{1}{\tau}\int_{E\triangle F}d_{\varphi^\circ_{hex}}^\epsilon(x,\partial F)\,\mathrm{d}x\\
&= P_\epsilon(E)+\frac{\sqrt{3}}{2}\frac{\epsilon^2}{\tau}\left(\sum_{H_\epsilon(i)\in E\backslash F} d_{\varphi^\circ_{hex}}^\epsilon(H_\epsilon(i),F)+\sum_{H_\epsilon(i)\in F\backslash E} d_{\varphi^\circ_{hex}}^\epsilon(H_\epsilon(i),\mathbb{R}^2\backslash F)\right),
\end{split}
\label{newenergy}
\end{equation}
where $\frac{\sqrt{3}}{2}\epsilon^2$ is the area of the hexagonal cell.

More precisely, given an initial set $E_0^{\tau,\epsilon}\in\mathcal{D}_\varepsilon$ approximating, as $\varepsilon,\tau\to0$ in the Hausdorff sense, a sufficiently regular set $E_0$, we define recursively a sequence $E^{\tau,\epsilon}_k$ in $\D_\epsilon$ by requiring that
$E^{\epsilon,\tau}_{k+1}$ is a minimizer of the functional $\mathcal{F}_{\tau,\epsilon}(\cdot,E^{\tau,\epsilon}_k)$, for every $k\geq0$.
Then, setting 

\begin{equation}\label{disefo}
E^{\tau,\epsilon}(t)=E^{\tau,\epsilon}_{\lfloor t\slash\tau\rfloor},
\end{equation}
for every $t\geq0$, we are interested in characterizing the motion described by any converging subsequence of $E^{\tau,\epsilon}(t)$ as $\epsilon,\tau\to0$. 


As remarked in the Introduction, the interaction between the two discretization parameters, in time and space, plays a crucial role in such a limiting process. More precisely, the limit motion depends strongly on their relative decrease rate to 0. If $\tau/\epsilon\to+\infty$, then we may first let $\epsilon\to0$, so that $P_\epsilon(E)$ can be directly replaced by the limit anisotropic perimeter $P(E)$ defined in (\ref{functionalhexa}) and $\frac{1}{\tau}\int_{E\bigtriangleup F}d_{\varphi^\circ_{hex}}^\epsilon(x,\partial F)\,\mathrm{d}x$ by 

\begin{equation}
\frac{1}{\tau}\displaystyle\int_{E\bigtriangleup F}\inf_{y\in\partial F}\varphi^\circ_{hex}(x-y)\,\mathrm{d}x.
\end{equation}
\\
As a consequence, the approximated flat motions tend to the solution of the time-continuous ones studied by {Almgren} and {Taylor}~\cite{AT95}, Taylor~\cite{Ta} with natural mobility function $M=\varphi_{hex}$.

On the other hand, if $\tau/\epsilon\to0$ then there is no motion (`pinning') since $E^{\tau,\epsilon}_k\equiv E_0^{\tau,\e}$. Indeed, for any $F\neq E_0^{\tau,\e}$ and for $\tau$ small enough we have

\begin{equation*}
\frac{1}{\tau}\int_{E_0^{\tau,\e}\bigtriangleup F}d_{\varphi^\circ_{hex}}^\epsilon(x,\partial F)\,\mathrm{d}x\geq C\frac{\epsilon}{\tau}>P_\epsilon(E_0^{\tau,\e}).
\end{equation*}

In this case, the limit motion is the constant state $E_0$. Hence, the meaningful regime is the intermediate case $\tau/\epsilon\to\gamma\in(0,+\infty)$ and we will focus on this case in the next Section.

\bigskip

\section{Motion of a convex ``Wulff-like'' set}\label{hexagonal}

\begin{figure}[htbp]
\centering
\def\svgwidth{100pt}
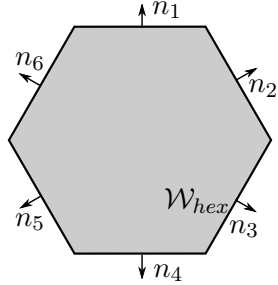
\caption{The Wulff shape $\mathcal{W}_{hex}$ of $\varphi_{hex}$ is a regular hexagon.}\label{fig:1}
\end{figure}

We introduce a class of sets, the convex ``Wulff-like'' sets, for which the motion by crystalline curvature {exists and is unique (at least) until the length of some side approaches to zero, and it is governed by a system of ordinary differential equations.} 
Roughly speaking, a convex {Wulff-like} set has a polygonal boundary that is a `good curve' made of regular corners, according to J. Taylor's terminology~\cite{Ta, AT95}; i.e., a convex set whose sides have the same exterior unit normal vectors and form the same angles as those of the Wulff shape $\mathcal{W}_{{hex}}$ of the density $\varphi_{hex}$. The Wulff shape $\mathcal{W}_{hex}$ is the regular hexagon 
\begin{equation}
\mathcal{W}_{hex}=\{\varphi_{hex}^\circ\leq1\}=\frac{4}{{3}}conv(\pm \eta_1,\pm \eta_2, \pm \eta_3),
\label{wulff}
\end{equation}
as pictured in Fig.~\ref{fig:1}.

To simplify the notation, we relabel in clockwise order the exterior unit normal vectors of the Wulff shape  $\mathcal{W}_{hex}$ (\ref{wulff}) and we set
\begin{equation}
\mathcal{N}:=\{n_1,n_2,n_3,n_4,n_5,n_6\},
\end{equation}
where $n_1:=\eta^\perp_1$, $n_2:=\eta^\perp_2$, $n_3:=\eta^\perp_3$, $n_4:=-n_1$, $n_5:=-n_2$, $n_6:=-n_3$.

Moreover, we will denote simply by $d^\epsilon$ the discrete distance $d_{\varphi^\circ_{hex}}^\epsilon$ defined in (\ref{discretedistance}).

\begin{defn}[Wulff-like set]\label{wshlikeset}
A bounded set $E\subset\mathbb{R}^2$ is said to be \emph{Wulff-like} if its boundary $\partial E$ is a polygonal closed curve whose sides $S_i, i=1,\dots,m$ have exterior unit normal vectors $\nu_i$ such that
\begin{enumerate}
\item[(1)] $\nu_i\in\mathcal{N}$, for every $i=1,\dots,m$;
\item[(2)] if $\nu_i=n_j$ for some $j=1,\dots,6$, then $\nu_{i+1}\in\{n_{j-1}, n_{j+1}\}$.
\end{enumerate}
Here, the labellings of $\nu_i$ and $n_j$ are intended to be modulo $m$ and $6$, respectively.
\end{defn}

\begin{figure}[htbp]
\centering
\def\svgwidth{150pt}
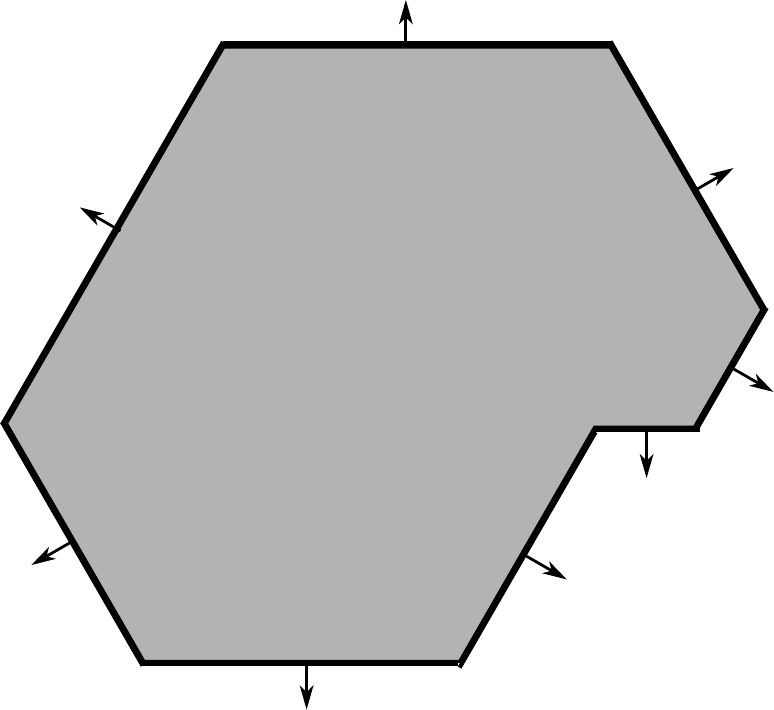
\caption{An example of Wulff-like set.}\label{fig:2}
\end{figure}

{Equivalently, each pair of adjacent sides of a Wulff-like set forms either a \emph{regular} or an \emph{inverse corner}, according to the definition given in \cite[Section~4]{AT95}. In particular, each side is parallel to a side of the Wulff shape. } 
Moreover, as an easy remark we note that \emph{convex} Wulff-like sets are convex Wulff-like hexagons. Among these, origin-symmetric convex Wulff-like hexagons will play the same role of rectangles on the square lattice in \cite{BGN}, and we will see that this case contains the main features of the motion.
Note that the evolution of more general (non-convex) sets may be studied up to assign a curvature sign on each side (see \cite[Sections 3.2-3.3]{BGN}).

{We will restrict the minimization of the energy $\mathcal{F}_{\tau,\varepsilon}$ defined in \eqref{newenergy} to those sets in $\mathcal{D}_\epsilon$ that are the union of all the cells of the hexagonal lattice strictly contained in a given convex Wulff-like hexagon. With a slight abuse of notation, we call such sets \emph{discrete} convex Wulff-like hexagons.}

{\begin{defn}[discrete convex Wulff-like hexagon]\label{dcwsl}
Let $\mathcal{D}_\varepsilon$ be defined as in \eqref{admissibleset}. A set $E\in\mathcal{D}_\epsilon$ is said to be a \emph{discrete convex Wulff-like hexagon} if there exists a convex Wulff-like hexagon $K$ such that
\begin{equation}
E=\bigcup_i\bigl\{H_\epsilon(i):\, H_\epsilon(i)\subset K\bigr\}.
\end{equation}
We denote this subclass by $\widetilde{\mathcal{D}}_\epsilon$.
\end{defn}}

{\begin{defn}[Wulff-like envelope]
Given any $E\in\widetilde{\mathcal{D}}_\epsilon$, we define $\mathcal{W}(E)$ the \emph{Wulff-like envelope} of $E$ as the smallest convex Wulff-like hexagon containing $E$. 
\end{defn}}

\begin{figure}[htbp]
\centering
\def\svgwidth{150pt}
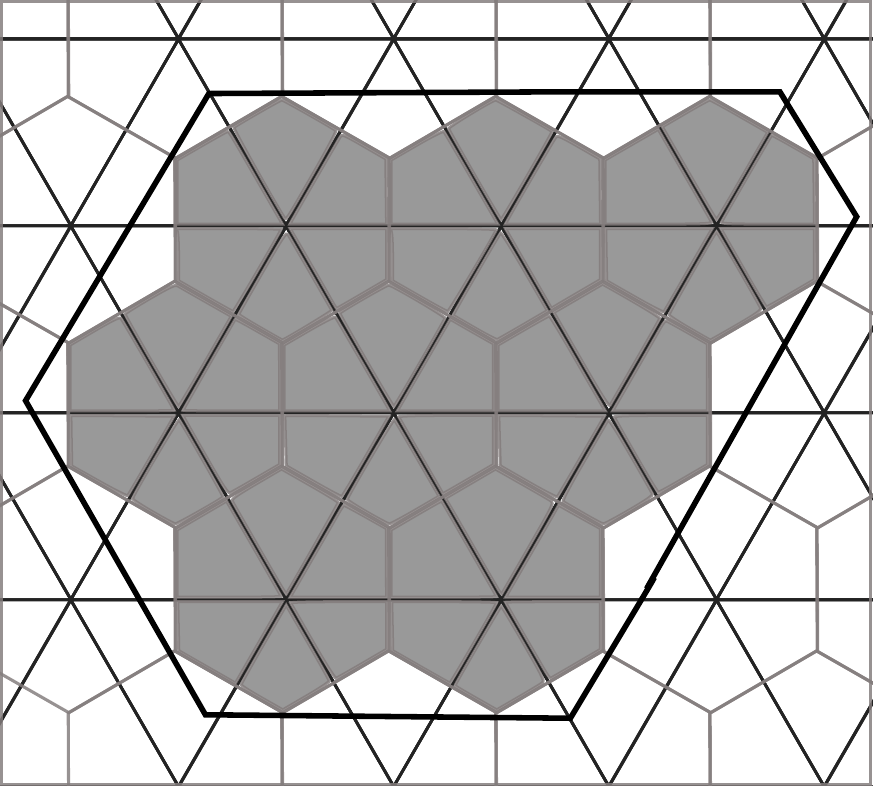
\caption{An example of discrete convex Wulff-like set.}\label{fig:discrete}
\end{figure}

An example of discrete convex Wulff-like hexagon is pictured in Fig.~\ref{fig:discrete}. The set is coloured in gray, while the continuous black line represents the boundary of its Wulff-like envelope.\\ 

Let $\tau=\tau(\epsilon)$. We characterize the microscopic motion of discrete convex Wulff-like hexagons at the critical regime
\begin{equation}
\lim_{\varepsilon\to0}\frac{\tau}{\varepsilon}=\lambda\in(0,+\infty).
\label{critreg}
\end{equation}
It is not restrictive to assume $\tau=\gamma\epsilon$ in place of \eqref{critreg}, since this only requires minor changes in the proof. Correspondingly, we omit the dependence on $\tau$ in the notation of $E_k^{\tau,\epsilon}=E_k^{\gamma\epsilon,\epsilon}$, that will be simply denoted by $E_k^\epsilon$.

The main result is that discrete convex Wulff-like hexagons evolve into sets of the same type. That is the content of Proposition~\ref{rectangleprop}, whose proof is reminiscent of some geometric arguments developed in \cite[Theorem~1]{BGN} for rectangular evolutions with underlying square lattice. However, we have to face some technical difficulties due to the geometry of the triangular lattice. Namely, (i) the boundary of a set in $\widetilde{\mathcal{D}}_\epsilon$ and that of its Wulff-like envelope do not coincide; (ii) in general, a convex Wulff-like hexagon is not axially symmetric and its center cannot be defined.

\begin{prop}\label{rectangleprop}
If $E_0^{\epsilon}\in{\widetilde{\D}}_\epsilon$ is a discrete convex Wulff-like hexagon and $E_k^\epsilon$ is a minimizer for the minimum problem for $\mathcal{F}_{\tau,\epsilon}(\cdot,E_{k-1}^\epsilon), k\geq1$, then 
$E_k^\epsilon$ is a discrete convex Wulff-like hexagon contained in $E_{k-1}^\epsilon$ as long as there exists $\delta>0$ such that the sides of its Wulff-like envelope $\mathcal{W}(E_{k-1}^\epsilon)$ are larger than $\delta$. 
\end{prop}

\begin{figure}[ht]
\centering
\def\svgwidth{150pt}
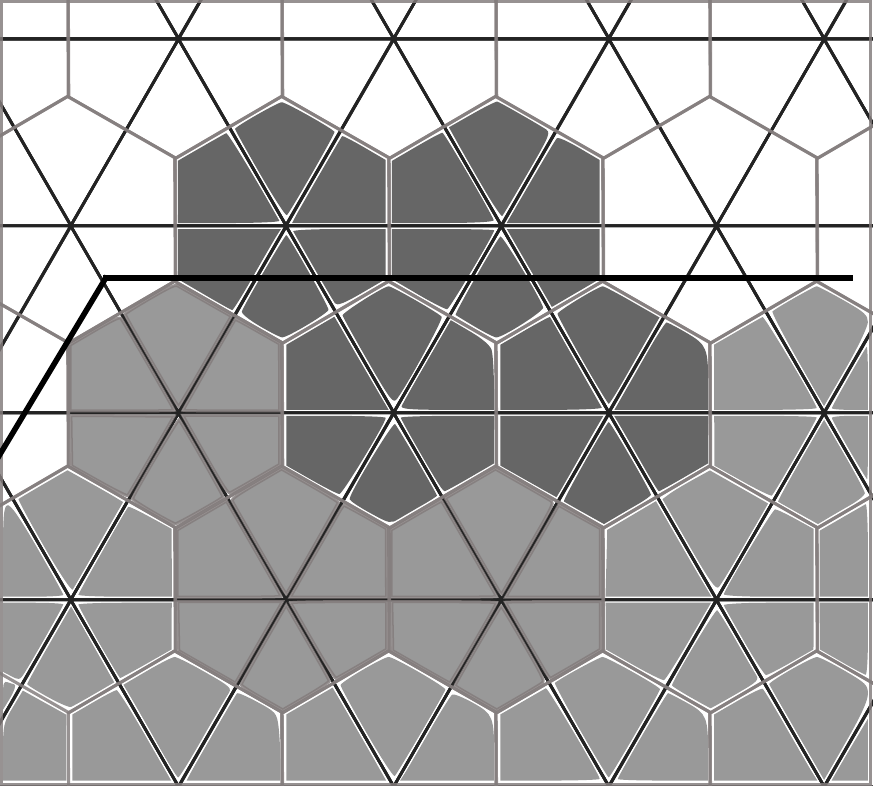\hfil
\def\svgwidth{150pt}
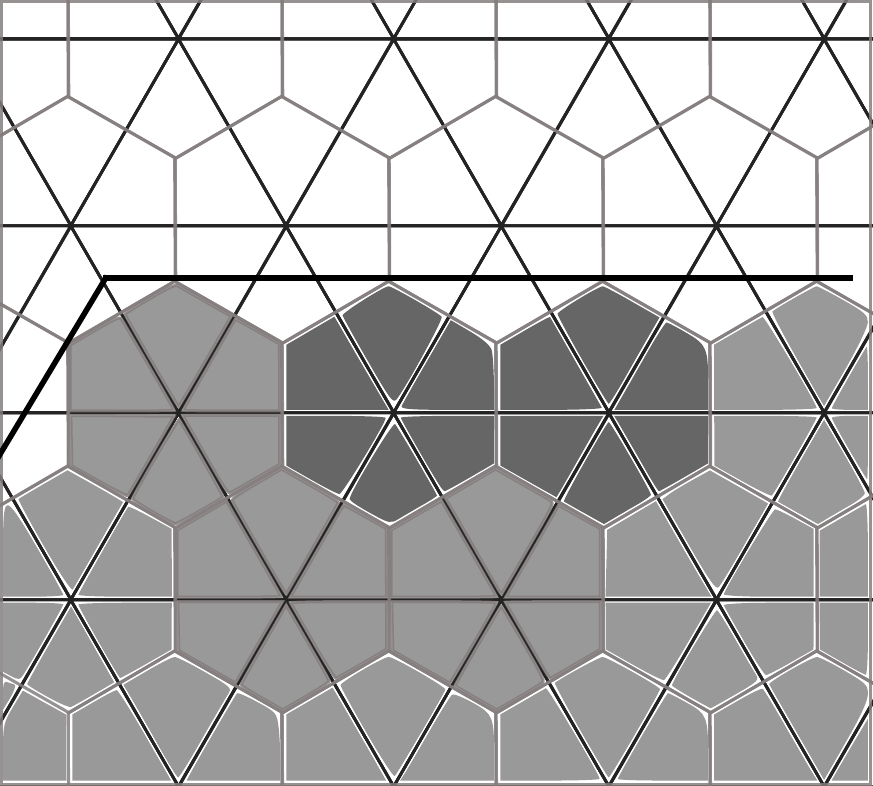
\caption{Each $F_i$ is contained in $E^\epsilon_0$.}\label{contain}
\end{figure}

\proof
The existence of minimizers among the sets of finite perimeter relies on classical results of compactness and semicontinuity (see, e.g., \cite[Section~3.2]{ATW83}). Here we characterize the geometrical properties of a minimizer. For this, it will suffice to show the assertion for $F=E^\epsilon_1$ a minimizer of $\mathcal{F}_{\tau,\epsilon}(\cdot,E_{0}^\epsilon)$, since the general case will follow by induction on the step $k$. In order to do that, let $F=F_1\cup\dots\cup F_m$ be the decomposition of $F$ into its connected components.
\\
{\bf Step 1: each $F_i$ is a discrete convex Wulff-like hexagon contained in $E^\epsilon_0$.} First, we characterize each connected component of $E^\epsilon_1$. We note that $E^\epsilon_1\subseteq E^\epsilon_0$. If not, let $F_i$ be a connected component of $E^\epsilon_1$ such that $F_i\cap(E^\epsilon_0)^c\neq\emptyset$; if $F_i\subseteq (E^\epsilon_0)^c$, then we may strictly reduce the energy $\mathcal{F}_{\tau,\epsilon}(\cdot, E^\epsilon_0)$ simply by dropping it. If not, we could consider as a competitor the set $F_i\cap E^\epsilon_0$: the area clearly decreases and the same holds for the perimeter, since any external connected curve made by cell sides connecting any two points of $\partial F_i\cap\partial \mathcal{W}(E^\epsilon_0)$ and containing cell sides orthogonal to $\partial\mathcal{W}(E^\epsilon_0)$ has perimeter not smaller than the one determined by the path along $\partial E^\epsilon_0$ (see Fig.~\ref{contain}).

\begin{figure}[ht]
\centering
\def\svgwidth{180pt}
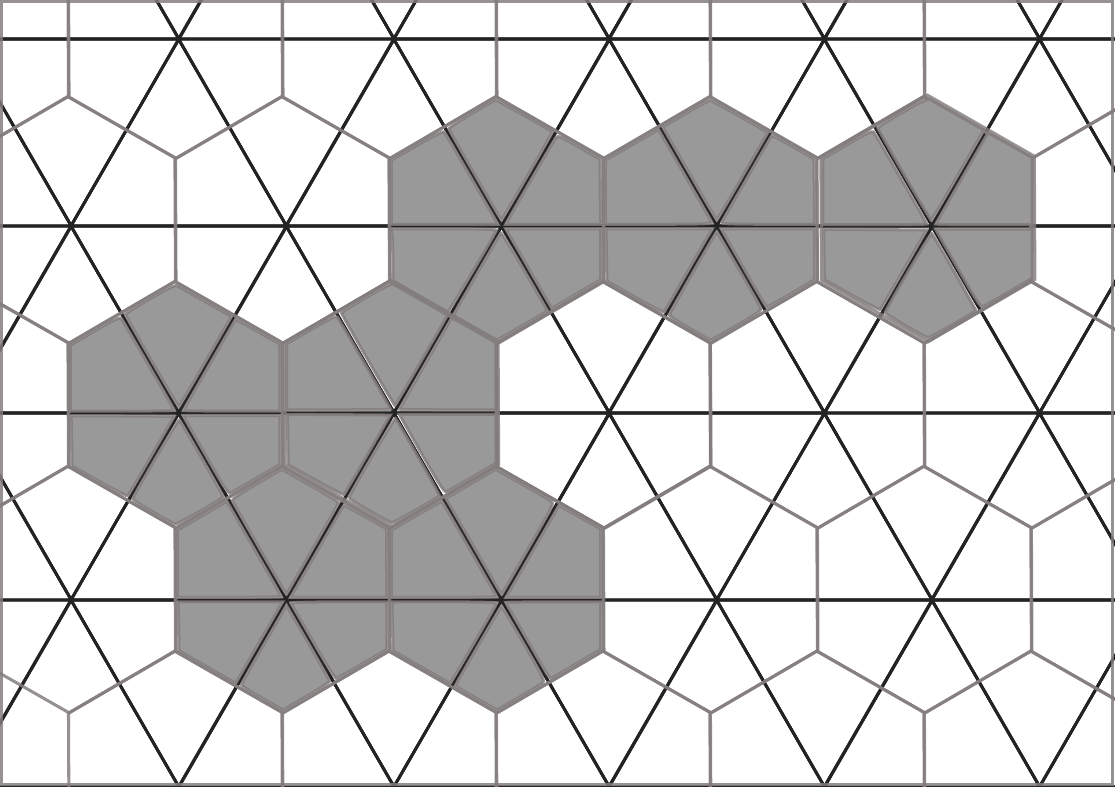\hfil
\def\svgwidth{180pt}
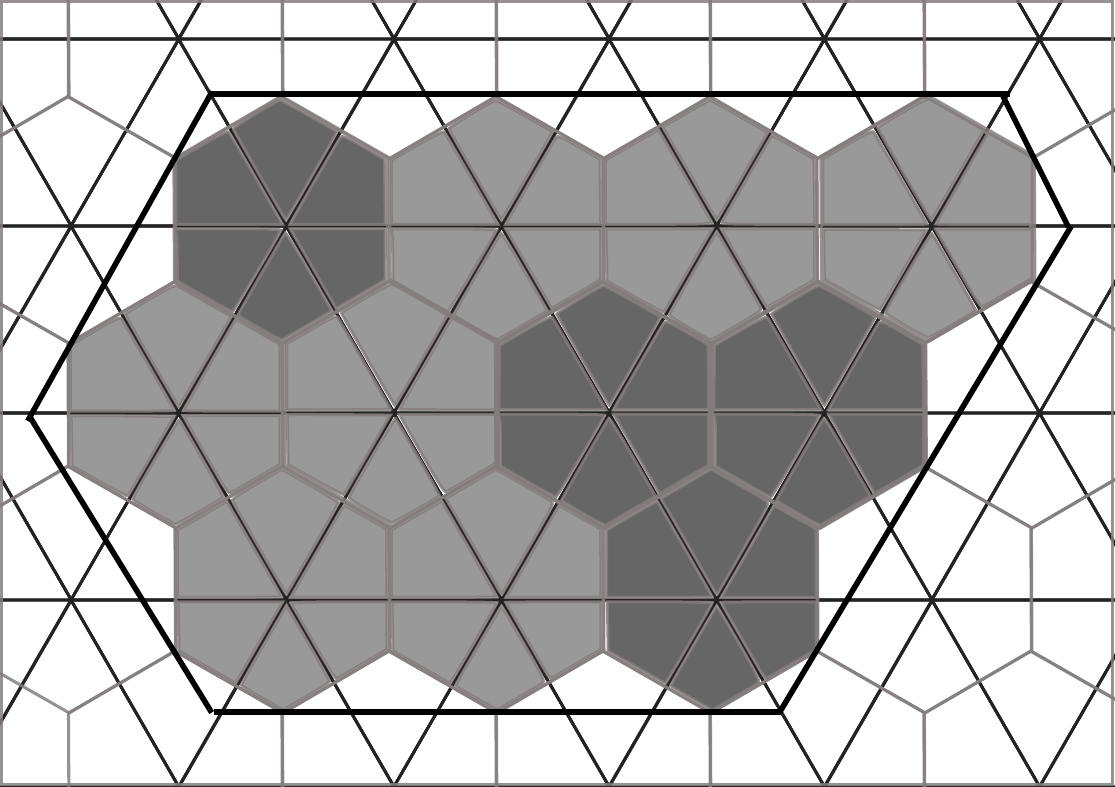
\caption{The Wulff-like ``convexification'' of $F_i$ reduces both its perimeter and the symmetric difference with the previous set.}\label{convex}
\end{figure}

Now, if $F_i$ is a discrete convex Wulff-like hexagon, then we are done. If not, we replace each $F_i$ with the smallest discrete convex Wulff-like hexagon containing $F_i$; in this case, its energy decreases since its perimeter is not greater than that of $F_i$ and the symmetric difference with $E^\epsilon_0$ decreases as well (see Fig.~\ref{convex}). More precisely, if we define

\begin{equation*}
\widetilde{F}_i=\cup_j\{H_\epsilon(j):\,H_\epsilon(j)\subseteq\mathcal{W}(F_i)\},
\end{equation*}
since $F_i\subseteq\widetilde{F}_i$ we immediately get that $|\widetilde{F}_i\triangle E^\epsilon_0|\leq|{F}_i\triangle E^\epsilon_0|$. In addition, we claim that $P_\epsilon(\widetilde{F}_i)\leq P_\epsilon({F}_i)$. For this, we consider the intersection of $\partial \widetilde{F}_i$ with a side $S_j$ of $\mathcal{W}(F_i)$ having normal vector $\pm\widetilde{\eta}_j$, for some $j=1,2,3$. Such an intersection will contain at least one point belonging to $\partial F_i\cap\partial \widetilde{F}_i$. Each connected curve made by cell sides passing through such points and different from the corresponding path along $\partial\widetilde{F}_i$ would contain some cell side (orthogonal to $S_j$) with normal vector $\pm \eta_j$, thus increasing the perimeter. Hence, if any portion of $\partial F_i$ is not contained in $\partial \widetilde{F}_i$, then $P_\epsilon({F}_i)\geq P_\epsilon(\widetilde{F}_i)$.
\\
\noindent
{\bf Step 2: each $F_i$ can be translated towards an incenter of $E^\epsilon_0$ without increasing its energy.} We rewrite the bulk term in the energy (\ref{newenergy}) as

\begin{equation*}
\begin{split}
\int_{F\bigtriangleup E_0^\e}d^\epsilon(x,\partial E_0^\e)\,\mathrm{d}x&=\int_{E_0^\e\backslash F}d^\epsilon(x,\partial E_0^\e)\,\mathrm{d}x=\int_{E_0^\e}d^\epsilon(x,\partial E_0^\e)\,\mathrm{d}x-\int_{F}d^\epsilon(x,\partial E_0^\e)\,\mathrm{d}x\\
&=\int_{E_0^\e}d^\epsilon(x,\partial E_0^\e)\,\mathrm{d}x-\int_{\cup_{i=1}^mF_i}d^\epsilon(x,\partial E_0^\e)\,\mathrm{d}x\\
&=\int_{E_0^\e}d^\epsilon(x,\partial E_0^\e)\,\mathrm{d}x-\sum_{i=1}^m\int_{F_i}d^\epsilon(x,\partial E_0^\e)\,\mathrm{d}x,
\end{split}
\end{equation*}
where in the latter equality the first integral is a fixed quantity independent of $F$ and the second integral is more negative if we translate along the coordinate directions each $F_i$ towards those points of $E^\epsilon_0$ with maximal distance from $\partial E^\epsilon_0$ (i.e., the kernel $\mathcal{M}(E^\epsilon_0)$ defined by (\ref{kernel})). As remarked in Section~\ref{setting}, $\mathcal{M}(E^\epsilon_0)$ is either a point (this is the case, for instance, when $\mathcal{W}(E^\epsilon_0)$ is a polygon with axial symmetry) or a line segment parallel to any of the sides of $\mathcal{W}(E^\epsilon_0)$.

We then distinguish between two cases:\\
{\bf (a) $\mathcal{M}(E^\epsilon_0)=\{C_0\}$.} We consider the component $F_1$, $C_0$ the incenter of $E^\epsilon_0$ and take a point $P\in F_1$. We can suppose that $P\neq C_0$, since there is at most one connected component of $F$ containing $C_0$. We define the set $F'$ obtained by substituting to $F_1$ its translation towards $C_0$ 
\begin{equation}
F'_1=F_1-\epsilon\,\mbox{sgn}(\langle P-C_0,\eta_1\rangle)\eta_1-\epsilon\,\mbox{sgn}(\langle P-C_0,\eta_2\rangle)\eta_2.
\label{translation1}
\end{equation}
The perimeter of $F_1'$ is the same as that of $F_1$, hence the perimeter term of $\mathcal{F}_{\tau,\epsilon}(F',E^\epsilon_0)$ remains unchanged, unless the boundary of $F'_1$ intersects the boundary of some other $F_j$ for a positive length (in which case the energy strictly decreases). We claim that the contribution of the bulk term in the energy does not increase, i.e.,
\begin{equation}
-\int_{F'_1}d^\epsilon(x,\partial E_0^\e)\,\mathrm{d}x\leq -\int_{F_1}d^\epsilon(x,\partial E_0^\e)\,\mathrm{d}x,
\label{estimate}
\end{equation}
by showing that if $P'=P-\epsilon\,\mbox{sgn}(\langle P-C_0,\eta_1\rangle)\eta_1-\epsilon\,\mbox{sgn}(\langle P-C_0,\eta_2\rangle)\eta_2$, then we have 
\begin{equation*}
d^\epsilon(P',\partial E_0^\e)\geq d^\epsilon(P,\partial E_0^\e).
\end{equation*}

\begin{figure}[htbp]
\centering
\def\svgwidth{200pt}
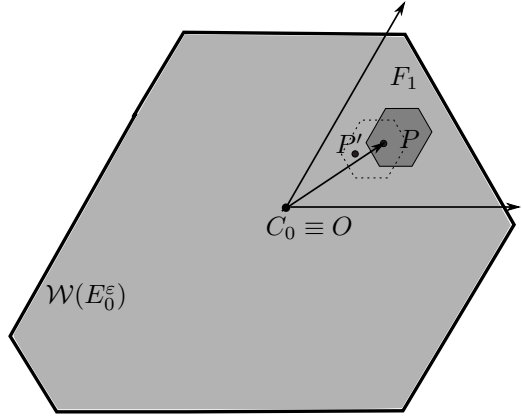
\caption{Each translation of $F_1$ towards $C_0$ does not increase the bulk term of the energy.}\label{figsystem}
\end{figure}

We take $C_0$ as the origin of a reference coordinate system (see Fig.~\ref{figsystem}) and, by the definition of $d^\epsilon$ as the piecewise constant interpolation of the corresponding values on lattice points, without loss of generality we can assume that $P\in\epsilon\mathbb{T}$, $P=(n\eta_1+m\eta_2)\epsilon$, $n,m\in\Z$. We can also assume that $n,m\geq0$, since the general case can be treated similarly. In order to simplify the computation, we view the translation defined by (\ref{translation1}) as the composition of a translation along $\eta_2$ ($n=0$) and an horizontal translation ($m=0$), showing that each of such elementary translations reduces (or leaves unchanged) the bulk term in the energy.

We first consider an elementary translation along $\eta_2$. In this case, $P=m\eta_2\epsilon$ and the shifted point is $P'=(m-1)\eta_2\epsilon$. If $P'\equiv C_0$ then trivially $d^\epsilon(P',\partial E_0^\e)\geq d^\epsilon(P,\partial E_0^\e)$ since, by definition, $C_0$ has maximal distance from $\partial E^\epsilon_0$. If not, let $q\in\N$ be such that $d^\epsilon(P,\partial E_0^\e)=\frac{3}{4}q\epsilon$; in this case $d^\epsilon(P',\partial E_0^\e)\in\{\frac{3}{4}q\epsilon,\frac{3}{4}(q+1)\epsilon\}\geq d^\epsilon(P,\partial E_0^\e)$. The same argument applies for elementary horizontal translations. Finally, the estimate (\ref{estimate}) follows by the arbitrariness of $P$.\\
\noindent
{\bf (b) $\mathcal{M}(E^\epsilon_0)$ is a line segment.} In this case, we can perform the same argument as in {\bf (a)} by choosing arbitrarily $C_0\in\mathcal{M}(E^\epsilon_0)$. 
\\
\noindent
{\bf Step 3: $F$ is connected.} If $m>1$ then the process at Step 2, applied to $F_1$ and $F_2$, after a finite number of steps produces a competitor $F'$ where the boundary of two such translated connected components, say $F'_1$ and $F'_2$, touch. Then their boundaries intersect in a set of positive length (a cell side), in which case a cancelation gives a lower contribution of the perimeter in contraddiction with the minimality of $F$. Hence, any minimizer $F$ has only one connected component, which is a discrete convex Wulff-like hexagon.

The previous construction can be iterated recursively for $k>1$, as long as $\mathcal{W}(E_{k-1}^\epsilon)$ remains an hexagon; that is, the length of each side of $\mathcal{W}(E_{k-1}^\epsilon)$ is greater than a positive constant $\delta$.
\endproof

In the following sections, we will compute explicitly the minimizer $E^\epsilon_k, k\geq1$ by a recursive minimization procedure that, in view of the latter result, can be performed among discrete convex Wulff-like hexagons contained in $E^\epsilon_{k-1}$.

\subsection{The pinning threshold}\label{pinning}

We first examine the case when the limit motion is trivial; i.e., all $E^\e_k$ are the same after a finite number of steps.
In case of rectangular evolutions with underlying lattice $\varepsilon\mathbb{Z}^2$ (see \cite{BGN, BraSci}), this is done by computing the {\em pinning threshold}; i.e., the critical value of the side length $L$ above which it is energetically not favorable for a side to move. This is obtained by imposing that the minimal displacement of a side by $\epsilon$ along any of the coordinate directions $e_1, e_2$ gives a non-negative contribution in the energy $\mathcal{F}_{\tau,\varepsilon}$. 

Also in our setting the coordinate directions of the underlying lattice coincide with the `preferred' directions for the motion $\{n_i,\,i=1,\dots,6\}$. 
The microscopic motion of a side is obtained by overcoming energy barriers along its normal direction; thus, the pinning threshold can be defined as \emph{the critical value $\overline{L}$ for the length of a side of the initial limit set above which it is not energetically favorable for such a side to move}.\\
\begin{figure}[htbp]
\centering
\def\svgwidth{190pt}
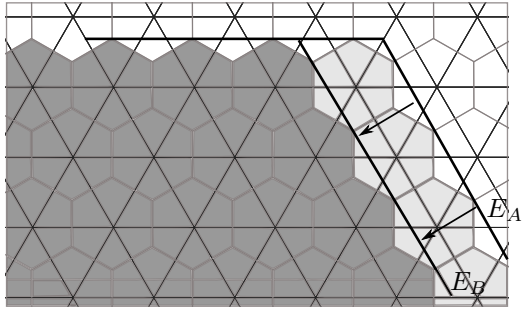
\caption{The minimal displacement of a side along its normal direction.}
\label{fig:pinning}
\end{figure} 

In order to determine it, we write the variation of the energy functional $\mathcal{F}_{\tau,\epsilon}$ from configuration $E_A$ to configuration $E_B$ in Fig.~\ref{fig:pinning}, regarding the inward translation by $\frac{\sqrt{3}}{2}\epsilon$ of a side of length $L$ along its normal direction; this consists in computing the variation of the energy corresponding to the removal of the layer of hexagonal cells coloured in light gray. If we impose it to be positive, we get
\begin{equation*}
\begin{split}
P_\epsilon(E_B)-P_\epsilon(E_A)+\frac{{3}}{4\tau}|E_A\backslash E_B|\epsilon&=-\frac{2\sqrt{3}}{3}\epsilon+\frac{{3}}{4\gamma\epsilon}\Bigl(\frac{L}{\epsilon}\frac{\sqrt 3}{2}\epsilon^2\Bigr)\epsilon\\
&=\sqrt{3}\epsilon\left[-\frac{2}{3}+\frac{{3}}{8}\frac{L}{\gamma}\right]\geq0,
\end{split}
\end{equation*}
whence we deduce that
\begin{equation}
L\geq\overline{L}=\alpha_{hex}\gamma:=\frac{16}{9}\gamma.
\label{threshold}
\end{equation}

\subsection{Description of the limit motion.}\label{limitmotion}

In this section, we provide the characterization of any limit motion in the critical regime $\tau=\gamma\epsilon$ for the class of \emph{origin-symmetric} convex Wulff-like hexagons. The restriction to this class ensures more symmetry for the motion. Indeed, with the following Theorem~\ref{limitmotion1}, we prove the convergence of the discrete scheme \eqref{newenergy}, as $\varepsilon\to0$, to a limit evolution which is a set of the same type. In particular, a regular Wulff-like hexagon shrinks homothetically until the extinction time. Among the main features of the limit evolutions due to the discrete motion described in Section~\ref{pinning}, we mention the phenomenon of ``quantization'' of the velocities. In other words, we expect velocities depending in a discontinuous way on the curvature. 

\begin{thm}\label{limitmotion1}For all $\epsilon>0$, let $E_0^\epsilon\in\widetilde{\D}_\epsilon$ be  discrete origin-symmetric convex Wulff-like hexagons and let the corresponding Wulff-like envelopes $\mathcal{W}(E_0^\epsilon)$ have sides $S^0_{1,\epsilon},\dots,S^0_{6,\epsilon}$. Assume also that
\begin{equation}
{\rm{d}}_\H(\mathcal{W}(E_0^\epsilon),E_0)<\epsilon
\label{approxinitial}
\end{equation}
for some fixed origin-symmetric convex Wulff-like hexagon $E_0$. Let $\gamma>0$ be fixed, 
$E^\e(t):=\mathcal{W}(E^\epsilon_{\lfloor t/\gamma\epsilon\rfloor})$ be the Wulff-like envelope of $E^\epsilon_{\lfloor t/\gamma\epsilon\rfloor}$ the piecewise-constant motion with initial datum $E_0^\e$ defined in {\rm(\ref{disefo})}. Then there exists $T>0$ such that $E^\epsilon(t)$ converges, up to a subsequence, as $\epsilon\to0$, in the Hausdorff topology and locally uniformly on $[0,T)$, to an origin-symmetric convex Wulff-like hexagon $E(t)$ with sides $S_i(t), i=1,\dots,6$ and such that $E(0)=E_0$. The distance $s_i(t)$ of the side $S_i(t)$ from the origin $O$ reduces with a velocity $v_i(t)$ satisfying
\begin{equation}
v_{i}(t)
\begin{cases}
=\displaystyle\frac{\sqrt3}{2}\Bigl(\frac{1}{\gamma}\left\lfloor\frac{\alpha_{hex}\gamma}{L_i(t)}\right\rfloor\Bigr), & \mbox{if $\displaystyle\frac{\alpha_{hex}\gamma}{L_i(t)}\not\in\mathbb{N}$}\\
\\
\in\displaystyle\frac{\sqrt{3}}{2}\left[\displaystyle\frac{1}{\gamma}\left(\frac{\alpha_{hex}\gamma}{L_i(t)}-1\right),\displaystyle\frac{\alpha_{hex}}{L_i(t)}\right], & \mbox{if $\displaystyle\frac{\alpha_{hex}\gamma}{L_i(t)}\in\mathbb{N}$}.
\end{cases}
\label{vl}
\end{equation}
where $L_i(t):=\H^1(S_i(t))$ denotes the length of the side $S_i(t)$. 
Accordingly, the law for each $L_i(t)$ is
\begin{equation}\label{system3biss}
\dot{L}_i(t)=\frac{2}{\sqrt{3}}\left[v_i(t)-\left(v_{i-1}(t)+v_{i+1}(t)\right)\right],\quad i=1,\dots,6,
\end{equation}
with the convention that symbols $v_i$ and $v_j$ coincide if $i\equiv j$ modulo 6.
\end{thm}

Before entering into details of the proof, we note that a natural choice for $T>0$ (see \cite[Theorem~3]{BGN}) is the first time for which $\lim_{t\to T^-}L_i(t)=0$ for some $i\in\{1,2,\dots,6\}$. If from one hand it is worth mentioning that an extension of the crystalline motion for times past such $T$, consisting of both deleting vanishing sides and possibly merging some other sides, has been provided in \cite[Section~2.3]{Ta}, from the other hand we note that an analogous delicate construction in the discrete setting would be out of this paper's scope. Thus, we agree that the evolution of convex Wulff-like hexagons exists until one of its sides vanishes.\\

\noindent
{\bf\emph {Proof of Theorem~\ref{limitmotion1}}.} 
{In view of Proposition~\ref{rectangleprop}, each $E^\epsilon_k$, $k\geq1$, is a discrete convex Wulff-like hexagon contained in $E^\epsilon_{k-1}$. We provide here the explicit computation of the minimizer $E^\epsilon_1$, since it can be iterated recursively at each step $k>1$.}  First, we note that if $E^\epsilon_0$ is origin-symmetric, then $E^\epsilon_1$ contains the origin $O$. If not, as in the proof of \cite[Theorem~1]{BGN}, we may consider the discrete convex Wulff-like hexagon $F$ being the symmetric of $E^\epsilon_1$ with respect to the origin. In this case, by symmetry, we have $P_\epsilon(F)=P_\epsilon(E^\epsilon_1)$ and $\int_{F\bigtriangleup E_0^\e}d^\epsilon(x,\partial E_0^\e)\,\mathrm{d}x=\int_{E^\epsilon_1\bigtriangleup E_0^\e}d^\epsilon(x,\partial E_0^\e)\,\mathrm{d}x$. Moreover, a comparison between the values $\mathcal{F}_{\tau,\epsilon}(E_1^\epsilon,E_{0}^\epsilon)$ and $\mathcal{F}_{\tau,\epsilon}(\emptyset,E_{0}^\epsilon)$ gives
\begin{equation*}
P_\epsilon(E^\epsilon_1)\leq\frac{1}{\tau}\int_{E^\epsilon_1}d^\epsilon(x,\partial E_0^\e)\,\mathrm{d}x,
\end{equation*}
whence
\begin{equation*}
\mathcal{F}_{\tau,\epsilon}(F\cup E_1^\epsilon,E_{0}^\epsilon)=\mathcal{F}_{\tau,\epsilon}(E_1^\epsilon,E_{0}^\epsilon)+P_\epsilon(E_1^\epsilon)-\frac{1}{\tau}\int_{F}d^\epsilon(x,\partial E_0^\e)\,\mathrm{d}x\leq \mathcal{F}_{\tau,\epsilon}(E_1^\epsilon,E_{0}^\epsilon).
\end{equation*}
Thus, $F\cup E_1^\epsilon$ is also a minimizer, and this contradicts the connectedness of the minimizer provided by Proposition~\ref{rectangleprop}.

Now, we prove that $E_1^\epsilon$ is origin-symmetric, by explicitly computing it. Let $L_{i}^\epsilon:=|S_{i}^{0,\epsilon}|$ be the length of the $i$-th side $S_{i}^{0,\epsilon}$ of $\mathcal{W}(E^\epsilon_0)$, $S_i^{1,\epsilon}$ be the $i$-th side of $\mathcal{W}(E^\epsilon_1)$ and $L_i^{1,\epsilon}$ its length, $\frac{\sqrt{3}}{2}\epsilon \overline{N}_i$, with $\overline{N}_i$ integer, and $s_i^\epsilon$ be the distances of the side $S_{i}^{1,\epsilon}$ from $S_i^{0,\epsilon}$ and $O$, respectively. If we subdivide the area between $S_{i}^{1,\epsilon}$ and $S_i^{0,\epsilon}$ in $\overline{N}_{i}$ layers of $L_i^\epsilon/\varepsilon$ hexagonal cells indexed by $k$, for each of which the discrete distance from the boundary is $\frac{3}{4}k\epsilon$, we can write the functional $\mathcal{F}_{\tau,\epsilon}(\cdot,E^\epsilon_0)$ in terms of the integers $N_{i}$, and we get that $\overline{N}_{1},\overline{N}_{2},\dots,\overline{N}_{6}$ are the minimizers of the function
\begin{equation}
\begin{split}
f(N_{1},N_{2},\dots,N_{6})&=-\frac{2\sqrt{3}}{3}\epsilon\sum_{i=1}^6 N_{i}+\frac{\epsilon^2}{\tau}\sum_{i=1}^6\sum_{k=1}^{N_{i}}\frac{{3}}{4}(k\epsilon)\frac{L_{i}^\epsilon}{\epsilon}\frac{\sqrt{3}}{2}-\frac{\epsilon^2}{\gamma}e_\epsilon\\
&=\sqrt{3}\epsilon\sum_{i=1}^6\left(-\frac{2}{3}N_{i}+\frac{3}{8\gamma}\frac{N_{i}(N_{i}+1)}{2}L_i^\epsilon\right)-\frac{\epsilon^2}{\gamma}e_\epsilon,
\end{split}
\label{functional}
\end{equation}
\\
with $0<e_\epsilon\leq C\max(\overline{N}_1,\dots,\overline{N}_6)^3$. The error $e_\epsilon$ is due to the bulk contribution of the hexagons near the vertices of $S_{i,\epsilon}$, which is negligible as $\epsilon\to0$.

The minimizer of \eqref{functional} is characterized by the inequalities
\begin{equation*}
f(\dots, \overline{N}_i,\dots)\leq f(\dots, \overline{N}_i\pm1,\dots),\quad i=1,\dots,6.
\end{equation*}
Since $f(\dots, N_i,\dots)$ is a parabola, the optimal value $\overline{N}_i$ of $N_i$ is the integer closest to
\begin{equation*}
\frac{\alpha_{hex}\gamma}{L_i^\epsilon}-\frac{1}{2},
\end{equation*} 
that is, $\overline{N}_i=\left\lfloor\frac{\alpha_{hex}\gamma}{L_i^\epsilon}\right\rfloor$ unless $\frac{\alpha_{hex}\gamma}{L_i^\epsilon}$ lies in a small neighborhood of the integers, infinitesimal as $\epsilon\to0$, when both an integer and the subsequent one are minimizers.

More precisely, as discussed in \cite[p. 480]{BGN}, there exists a constant $\bar{C}=\bar{C}(L_1,\dots,L_6)$ such that this problem has the unique minimizer $\overline{N}_i=\lfloor\frac{\alpha_{hex}\gamma}{L_i^\epsilon}\rfloor$ if $\dist\left(\frac{\alpha_{hex}\gamma}{L_i^\epsilon},\mathbb{N}\right)\geq \bar{C}\epsilon$, otherwise we have a double choice for $\overline{N}_i$. Namely,
\begin{equation*}
\begin{split}
\overline{N}_i\in\left\{\left\lfloor\frac{\alpha_{hex}\gamma}{L_i^\epsilon}\right\rfloor-1,\left\lfloor\frac{\alpha_{hex}\gamma}{L_i^\epsilon}\right\rfloor\right\}, &\quad{\mbox{ if }}\,\,\frac{\alpha_{hex}\gamma}{L_i^\epsilon}-\left\lfloor\frac{\alpha_{hex}\gamma}{L_i^\epsilon}\right\rfloor< \bar{C}\epsilon;\\
\overline{N}_i\in\left\{\left\lfloor\frac{\alpha_{hex}\gamma}{L_i^\epsilon}\right\rfloor,\left\lfloor\frac{\alpha_{hex}\gamma}{L_i^\epsilon}\right\rfloor+1\right\}, &\quad{\mbox{ if }}\,\,\left\lfloor\frac{\alpha_{hex}\gamma}{L_i^\epsilon}\right\rfloor+1-\frac{\alpha_{hex}\gamma}{L_i^\epsilon}< \bar{C}\epsilon.
\end{split}
\end{equation*}

This singularity affects the uniqueness of the limit velocity of each side. Indeed, if in the limit as $\epsilon\to0$, $\frac{\alpha_{hex}\gamma}{L_i}\not\in\mathbb{N}$, then the velocity of the $i$-th side is uniquely determined by $\frac{\sqrt{3}}{2\gamma}\left\lfloor\frac{\alpha_{hex}\gamma}{L_i}\right\rfloor$. If, instead, $\frac{\alpha_{hex}\gamma}{L_i}\in\mathbb{N}$, then the velocity is not unique since it depends how the value $\frac{\alpha_{hex}\gamma}{L_i^\epsilon}$ approaches $\mathbb{N}$ at the $\epsilon$ level. However, also in this case a limit velocity can be defined (see \cite{BGN} for details) leading to formulas \eqref{vl}.

We have
\begin{eqnarray*}
&s_i^{1,\epsilon}=s_i^{\epsilon}-\frac{\sqrt3}{2}\overline{N}_{i}\epsilon,\\
&L_i^{1,\epsilon}=L_i^{\epsilon}+\bigl(\overline{N}_{i}-(\overline{N}_{i-1}+\overline{N}_{i+1})\bigr)\epsilon, 
\end{eqnarray*}
where, by the symmetry assumption on $E_0^\epsilon$, $\overline{N}_{1}=\overline{N}_{4}$, $\overline{N}_{2}=\overline{N}_{5}$ and $\overline{N}_{3}=\overline{N}_{6}$, thus giving $L_1^{1,\epsilon}=L_4^{1,\epsilon}$, $L_2^{1,\epsilon}=L_5^{1,\epsilon}$ and $L_3^{1,\epsilon}=L_6^{1,\epsilon}$. In particular, $E_1^\epsilon$ is origin-symmetric.

The process described for $k=1$ can be iterated constructing recursively for $k>1$ and $i=1,\dots,6$ three sequences $s_i^{k,\epsilon}$, $L_i^{k,\epsilon}$ and $N_{i}^{k,\epsilon}$ such that

{\begin{eqnarray}
&s_i^{k+1,\epsilon}=s_i^{k,\epsilon}-\frac{\sqrt3}{2}N_{i}^{k,\epsilon}\epsilon,\label{incremental1}\\
&L_i^{k+1,\epsilon}=L_i^{k,\epsilon}+\frac{2}{\sqrt{3}}\bigl(s_{i+1}^{k+1,\epsilon}-s_{i+1}^{k,\epsilon}+s_{i-1}^{k+1,\epsilon}-s_{i-1}^{k,\epsilon}-(s_{i}^{k+1,\epsilon}-s_{i}^{k,\epsilon})\bigr), \label{incremental2}
\end{eqnarray}
\\
since by geometry there holds (see, e.g., \cite[p. 423]{Ta})
\begin{equation}
\begin{split}
L_i^{k,\epsilon}&=\frac{s_{i+1}^{k,\epsilon}-\langle n_{i+1}, n_i\rangle s_i^{k,\epsilon}}{\sqrt{1-(\langle n_{i+1}, n_i\rangle)^2}}+\frac{s_{i-1}^{k,\epsilon}-\langle n_{i-1}, n_i\rangle s_i^{k,\epsilon}}{\sqrt{1-(\langle n_{i-1}, n_i\rangle)^2}}\\
&=\frac{2}{\sqrt{3}}\bigl(s_{i+1}^{k,\epsilon}+s_{i-1}^{k,\epsilon}-s_{i}^{k,\epsilon}\bigr),
\end{split}
\label{geomconstr}
\end{equation}}
with initial conditions $s_i^{0,\epsilon}=s_i^\epsilon$, $N_{i}^{0,\epsilon}=\overline{N}_{i}$ and $L_i^{0,\epsilon}=L_i^\epsilon$. $N_{i}^{k,\epsilon}$ is a minimizer obtained by the same minimization procedure as above with $L_i^{k,\epsilon}$ in place of $L_i^{\epsilon}$. 

For each $1\leq i\leq6$, we define $\bar{s}_i^\epsilon(t)$ and $\bar{L}_i^\epsilon(t)$ the piecewise affine interpolations in $[k\tau, (k+1)\tau]$ of the values $s_i^{k,\epsilon}$ and $L_i^{k,\epsilon}$, respectively. From \eqref{incremental2} we deduce the identity
\begin{equation}
\bar{L}_i^\epsilon(t)=L_i^{0,\epsilon}+\frac{2}{\sqrt{3}}(s_i^{0,\epsilon}-\bar{s}_i^\epsilon(t)-(s_{i-1}^{0,\epsilon}-\bar{s}_{i-1}^\epsilon(t)+s_{i+1}^{0,\epsilon}-\bar{s}_{i+1}^\epsilon(t))).
\label{incremental3}
\end{equation} 

Let $T_\epsilon>0$ be defined as
\begin{equation}
T_\epsilon:=\sup\left\{t>0:\, \exists c>0\, \mbox{ such that } \bar{L}_i^\epsilon(r)\geq c,\quad\forall r\in[0,t),\, \mbox{ for every $i$}\right\},
\end{equation}
and let $\overline{T}\in(0,T_\epsilon)$ be arbitrarily fixed. By \eqref{incremental1} we have 

\begin{equation}
\frac{s_i^{k+1,\epsilon}-s_i^{k,\epsilon}}{\tau}=-\frac{\sqrt3}{2\gamma}N_{i}^{k,\epsilon},
\label{bound}
\end{equation}
\\
so that $\bar{s}_i^\epsilon(t)$ is a decreasing continuous function of $t$. Since by \eqref{approxinitial} we may assume that $|s_i^{0,\epsilon}- s_i^0|\leq c\epsilon$ for every $i$ and a suitable constant $c$, the monotonicity of $\bar{s}_i^\epsilon(t)$ implies the existence of $\epsilon_0>0$ and a uniform constant $C_1>0$ such that
\begin{equation}
|\bar{s}_i^\epsilon(t)|\leq C_1,\quad \mbox{ for every $t\in[0,\overline{T}]$, for $\epsilon\leq\epsilon_0$.}
\label{equibounded}
\end{equation} 
Moreover, it holds that
\begin{equation}
|\bar{s}_i^\epsilon(t_1)-\bar{s}_i^\epsilon(t_2)|\leq C_2 |t_1-t_2|,\quad \mbox{for every $t_1,t_2\in[0,\overline{T}]$,}
\label{equicontinuous}
\end{equation}
for some positive constant $C_2$ independent of $\epsilon$. In order to prove this, we may assume that $t_2<t_1$. Then, taking into account \eqref{bound}, we get
\begin{equation*}
\begin{split}
|\bar{s}_i^\epsilon(t_1)-\bar{s}_i^\epsilon(t_2)| & \leq |\bar{s}_i^\epsilon(t_1)-s_i^{\lfloor t_1/\tau\rfloor,\epsilon}|+\sum_{r=\lfloor t_2/\tau\rfloor+1}^{\lfloor t_1/\tau\rfloor-1}|s_i^{r+1,\epsilon}-s_i^{r,\epsilon}| + |s_i^{\lfloor t_2/\tau\rfloor+1,\epsilon}-\bar{s}_i^\epsilon(t_2)|\\
&\leq \frac{\sqrt3}{2\gamma} \left[(t_1-\lfloor t_1/\tau\rfloor\tau)+ (\lfloor t_1/\tau\rfloor-\lfloor t_2/\tau\rfloor-1)\tau+ ((\lfloor t_2/\tau\rfloor+1)\tau-t_2)\right]\\
&\leq \frac{\sqrt3}{2\gamma}|t_1-t_2|.
\end{split}
\end{equation*}
Hence, in view of \eqref{equibounded}-\eqref{equicontinuous}, by the Ascoli-Arzel\`a theorem there exists a subsequence $\bar{s}_i^{\epsilon_j}(t)$, with $\epsilon_j\to0$, converging uniformly on $[0,\overline{T}]$, as $j\to+\infty$, to a continuous function $s_i(t)$, which is also decreasing. 
Moreover, with \eqref{incremental3}, we get the convergence of $\bar{L}_i^\epsilon(t)$, as $\epsilon_j\to0$, 
to the function $L_i(t)$ defined as
\begin{equation*}
L_i(t):=L_i^0+\frac{2}{\sqrt{3}}(s_i^0-s_i(t)-(s_{i-1}^0-s_{i-1}(t)+s_{i+1}^0-s_{i+1}(t))),
\end{equation*}
where we used also the fact that, by \eqref{approxinitial} and \eqref{geomconstr},
\begin{equation*}
|L_i^{0,\epsilon}-L_i^0|\leq\frac{\sqrt{3}}{2}(|s_i^{0,\epsilon}-s_i^0|+|s_{i-1}^{0,\epsilon}-s_{i-1}^0|+|s_{i+1}^{0,\epsilon}-s_{i+1}^0|)\leq c'\epsilon.
\end{equation*}
Setting
\begin{equation}
T:=\sup\left\{t>0:\, {L}_i(r)>0,\quad\forall r\in[0,t),\, \mbox{ for every $i$}\right\},
\end{equation}
as a consequence of the convergence result, we have also that $\lim_{j}T_{\epsilon_j}=T$. This allows us to choose $\overline{T}$ arbitrarily close to the extinction time $T$.\\
It follows that $E^{\epsilon_j}(t)$ converges as $\epsilon_j\to0$, in the Hausdorff sense and locally uniformly on $[0,T)$, to the origin-symmetric convex Wulff-like hexagon $E(t)$ with sides of lengths $L_i(t)$, $i=1,2,\dots,6$, such that $E(0)=E_0$.\\

Now we justify the formula (\ref{vl}) for the velocities. To simplify the computation, we introduce the piecewise-constant interpolations of the values $s_i^{k,\epsilon},L_i^{k,\epsilon}, N_{i}^{\epsilon}$; namely, for $t\geq 0$ we put $s^\tau_i(t)=s_i^{\lfloor t/\tau\rfloor,\epsilon}$, $L^\tau_i(t)=L_i^{\lfloor t/\tau\rfloor,\epsilon}$ and $N_{i}^\tau(t)=N_{i}^{\lfloor t/\tau\rfloor,\epsilon}$. 
Note that, if $\epsilon,\tau\to0$, then
\begin{equation}
\bar{s}_i^\epsilon(t)-s^\tau_i(t)\to0,\quad \mbox{ uniformly with respect to }t\in[0,\overline{T}].
\end{equation}
Indeed, for every $t\in[k\tau,(k+1)\tau]$, from \eqref{bound} we have the estimate
\begin{equation}
|\bar{s}_i^\epsilon(t)-s^\tau_i(t)|=\frac{|s_i^{k+1,\epsilon}-s_i^{k,\epsilon}|}{\tau}|t-k\tau|\leq|s_i^{k+1,\epsilon}-s_i^{k,\epsilon}|\leq\frac{\sqrt{3}}{2}\epsilon.
\end{equation}
Thus, $s^\tau_i(t)\to s_i(t)$, $L^\tau_i(t)\to L_i(t)$ locally uniformly as $\tau\to0$ and, by continuity, $N_{i}^\tau(t)\to\frac{{2\gamma}v_{i}(t)}{\sqrt3}$ as $\tau\to0$, where the velocities $v_i(t)$ are defined by \eqref{vl}.
By construction we have
\begin{equation*}
\begin{split}
s_i^\tau(t+\tau)&=s_i^0-\frac{\sqrt3}{2\gamma}\sum_{k=0}^{\lfloor t/\tau\rfloor}\tau N_{i}^\tau(k\tau)\\
&=s_i^0-\sum_{k=0}^{\lfloor t/\tau\rfloor}\tau v_{i}(k\tau)+\omega(\tau),
\end{split}
\end{equation*}
$\omega(\tau)$ being an infinitesimal error as $\tau\to0$, where the second equality has been obtained using the convergence of $N_{i}^\tau$ to $\frac{{2\gamma}v_{i}(t)}{\sqrt3}$. Passing to the limit as $\tau\to0$ we finally deduce that
\begin{equation*}
s_i(t)=s_i^0-\int_0^tv_{i}(s)\,ds,
\end{equation*}
that is equivalent to (\ref{vl}) rephrased through the relation $\dot{s}_i(t)=-v_{i}(t)$.

As for \eqref{system3biss}, from \eqref{incremental2} and arguing as before we get
\begin{equation*}
\begin{split}
L_i^\tau(t+\tau)&=L_i^0+\sum_{k=0}^{\lfloor t/\tau\rfloor}\tau [N_{i}^\tau(k\tau)-(N_{i-1}^\tau(k\tau)+N_{i+1}^\tau(k\tau))]\\
&=L_i^0+\frac{2}{\sqrt{3}}\sum_{k=0}^{\lfloor t/\tau\rfloor}\tau [v_{i}(k\tau)-(v_{i-1}(k\tau)+v_{i+1}(k\tau))]+\omega(\tau),
\end{split}
\end{equation*}
whence, passing to the limit as $\tau\to0$, we obtain
\begin{equation*}
L_i(t)=L_i^0+\frac{2}{\sqrt{3}}\int_0^t v_i(s)-(v_{i-1}(s)+v_{i+1}(s))\,ds,
\end{equation*}
from which \eqref{system3biss} follows by taking the time derivative of both the sides.
\qed
\\

The following theorem characterizes the limit evolutions when we have a unique choice for the velocities $v_i$ in \eqref{vl}. That is the case, for instance, when (a) all the sides are pinned (``total pinning'') and the motion is trivial; (b) all the sides have the same length and are short enough (``motion of a Wulff shape''), thus obeying to a self-similar evolution that extinguishes in finite time.

\begin{thm}[unique limit motions]\label{unilimitmotion}
Let $E^\epsilon(t),E_0$ be as in the statement of Theorem~{\rm \ref{limitmotion1}}. Assume in addition that the lengths $L_i^0, i=1,\dots,6$ of the sides of the initial set $E_0$ satisfy one of the following conditions:
\begin{itemize}
\item[\emph{(a)}] $\displaystyle\min_{1\leq i\leq6}\{L_i^0\}>\alpha_{hex}\gamma$\,\, {\rm (total pinning);}
\item[\emph{(b)}] $L_i^0= L^0<\alpha_{hex}\gamma$ \,\,{\rm for\,\,every\,\,}$i=1,\dots,6$\,{\rm (self-similar evolution vanishing in finite time);} 
\end{itemize}
then there exists $T>0$ such that $E^\epsilon(t)$ converges locally in time to $E(t)$ on $[0,T)$ as $\epsilon\to 0$, where $E(t)$ is the unique origin-symmetric convex Wulff-like hexagon with sides $S_i(t)$ of length $L_i(t)$ whose distances from the origin $s_i(t)$, respectively, solve the following system of degenerate ordinary differential equation

\begin{equation}\label{unita}
\dot{s}_i(t)=-\frac{\sqrt{3}}{2}\left({1\over \gamma}\left\lfloor\frac{\alpha_{hex}\gamma}{L_{i}(t)}\right\rfloor\right), \quad i=1,\dots,6
\end{equation}
for almost every $t\in[0,T)$, with initial condition $L_i(0)=L_i^0$. Accordingly, the law for each side of length $L_i(t)$ is
\begin{equation}\label{system3bis}
\dot{L}_i(t)=\frac{1}{\gamma}\left\lfloor\frac{\alpha_{hex}\gamma}{L_{i}(t)}\right\rfloor-\frac{1}{\gamma}\left(\left\lfloor\frac{\alpha_{hex}\gamma}{L_{i-1}(t)}\right\rfloor+\left\lfloor\frac{\alpha_{hex}\gamma}{L_{i+1}(t)}\right\rfloor\right),\quad i=1,\dots,6,
\end{equation}
with the convention that symbols $L_i$ and $L_j$ coincide if $i\equiv j$ modulo 6.
\end{thm}

\proof
First, we note that \eqref{unita} follows by \eqref{vl}, while from \eqref{vl} and \eqref{system3biss} we deduce \eqref{system3bis}.\\
\noindent
(a) In this case the statement follows by \eqref{vl} noticing that we have $v_{i}(t)=0$ for all $t\geq0$, which is equivalent to $\dot{s}_i=0$. Correspondingly, $\dot{L}_i(t)\equiv0$. \\
(b) If we choose $L_i^0=L^0<\alpha_{hex}\gamma$ for every $i=1,\dots,6$, from \eqref{system3bis} we deduce that the side length $L(t)$ solves the differential equation with discontinuous right hand side
\begin{equation}\label{unitareg1}
\dot{L}(t)=-{1\over \gamma}\left\lfloor\frac{\alpha_{hex}\gamma}{L(t)}\right\rfloor,
\end{equation}
until the extinction time $T$. The law (\ref{unitareg1}) immediately follows also by \eqref{unita}, once we remark that $s_i(t)\equiv s(t)$ is independent of the sides and coincides with the apothem of $E(t)$, and it holds $\dot{L}(t)=\frac{2}{\sqrt{3}}\dot{s}(t)$. 
According to \eqref{unitareg1}, the side length $L(t)$ decreases, with a strictly negative derivative $\dot{L}(t)\leq-1/\gamma$, until it vanishes. Thus, $\frac{\alpha_{hex}\gamma}{L(t)}\in\mathbb{N}$ only for a countable set of times $t$ and uniqueness for the solution of equation \eqref{unitareg1} may be proved.

\endproof

The explicit description of the evolution \eqref{system3biss} for a very general convex Wulff-like hexagon is tricky, since the rate of change of the length for each side depends on the velocities of the neighboring sides. Furthermore, depending on the sign of the right hand side in \eqref{system3biss}, a side may shorten or lengthen, possibly reaching the pinning threshold after an initial motion. In this case, the uniqueness of velocities \eqref{vl} may no longer hold. Another feature we point out is the \emph{partial} pinning of a side; that is, the side stays pinned until it becomes sufficiently short, due to the motion of the adjacent sides, and then moves. We enlighten this phenomenon through the following example, where we consider a particular class of symmetric convex Wulff-like hexagons. 

\begin{figure}[htbp]
\centering
\def\svgwidth{300pt}
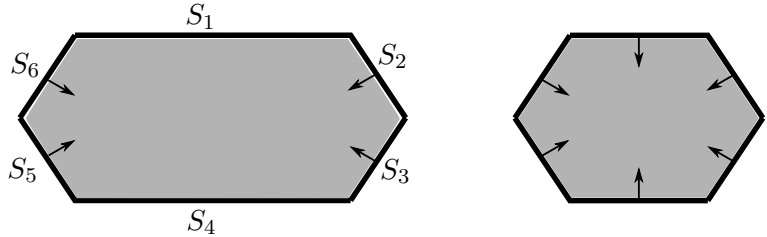
\caption{An example of partial pinning.}
\label{longset}
\end{figure}

\begin{oss}[An example of partial pinning] We consider a symmetric initial set where a pair of sides stays pinned, at least for a finite time (see Fig.~\ref{longset}). 

Let $L_1^0=L_{4}^0=:L^{0,1}>\alpha_{hex}\gamma$, $L_{6}^0=L_{2}^0=L_{3}^0=L_{4}^0=:L^{0,2}<\alpha_{hex}\gamma$, with $L^{0,2}<L^{0,1}$.
We examine the motion of sides $S_{1}, S_{2}, S_{3}$, since by symmetry the same holds for the triple $S_{4}, S_{5}, S_{6}$. We prove the following\\
\noindent
{\bf Claim:} There exists $T>0$ such that, for every $t\in[0,T]$, the side $S_1$ stays pinned, $L_{2}(t)=L_{3}(t)\equiv L^{0,2}$ and the length $L_1(t)$ reduces with constant velocity according to
\begin{equation}
\dot{L}_1(t)=-\frac{2}{\gamma}\left\lfloor\frac{\alpha_{hex}\gamma}{L^{0,2}}\right\rfloor,
\end{equation}
due to the motions of the neighboring sides $S_2$ and $S_6$. 

Indeed, at every time $t$ such that $v_1(t)=0$, the laws for the side lengths of $S_2$ and $S_3$ are
\begin{equation}
\dot{L}_2(t)=\frac{2}{\sqrt{3}\gamma}\bigl(v_2(t)-v_3(t)\bigr)
\label{eq1}
\end{equation}
and
\begin{equation}
\dot{L}_3(t)=\frac{2}{\sqrt{3}\gamma}\bigl(v_3(t)-v_2(t)\bigr),
\label{eq2}
\end{equation}
respectively, with initial conditions $L_2(0)=L_3(0)=L^{0,2}$. From \eqref{eq1}-\eqref{eq2} we deduce $\dot{L}_2(t)+\dot{L}_3(t)=0$, whence 
\begin{equation}
L_2(t)+L_3(t)=2L^{0,2}.
\label{eq3}
\end{equation}
Combining \eqref{eq3} with \eqref{eq1}, we obtain the differential equation
\begin{equation}
\dot{L}_2(t)=\frac{1}{\gamma}\left(\left\lfloor\frac{\alpha_{hex}\gamma}{L_2(t)}\right\rfloor-\left\lfloor\frac{\alpha_{hex}\gamma}{2L^{0,2}-L_2(t)}\right\rfloor\right)
\end{equation}
with initial condition $L_2(0)=L^{0,2}$, that admits the unique solution $L_2(t)\equiv L^{0,2}$. Thus, $L_2(t)=L_3(t)=L^{0,2}$, $v_2(t)=v_3(t)=\frac{\sqrt{3}}{2}\left\lfloor\frac{\alpha_{hex}\gamma}{L^{0,2}}\right\rfloor$ and $\dot{L}_1(t)=-\frac{4}{\sqrt{3}\gamma}v_2(t)$. Now, the time $T$ is determined by $L_1(T)=\alpha_{hex}\gamma$ and this concludes the proof of Claim. 

For $t>T$ the evolution is governed by equations
\begin{equation*}
\dot{L}_1(t)=\frac{2}{\sqrt{3}}\left(v_1(t)-2v_{2}(t)\right),
\label{lato1}
\end{equation*}
since $L_{2}(t)=L_{6}(t)$ for every $t$, where
\begin{equation*}
\dot{L}_{2}(t)=-\frac{2}{\sqrt{3}}v_1(t).
\end{equation*} 
\label{esempiopartial}
\end{oss}

The study of singular initial data, leading to some additional phenomena as \emph{non-uniqueness} of limit evolutions, would require a more refined analysis (we refer to \cite[Section~3.1.1]{BGN} for a discussion). We just mention that if the initial datum is a regular Wulff-like hexagon with side length $L_0=\alpha_{hex}\gamma$, then for every $T>0$ and up to choosing properly the discrete motions $E_k^\epsilon$, we can characterize the limit evolution as follows: for every $t\in[0,T]$, the initial hexagon stays pinned as in Theorem~\ref{unilimitmotion}(a); then, for every $t>T$, it shrinks homothetically to its center as in Theorem~\ref{unilimitmotion}(b).

\begin{oss}[a comparison with the crystalline motion with natural mobility]
According to Theorem~\ref{unilimitmotion}(b), if we choose $E_0$ to be a regular Wulff-like hexagon centered at the origin $O$, of side length $L_0$ complying with $L_0<\alpha_{hex}\gamma$, then the limit motion $E(t)$ is the unique regular Wulff-like hexagon centered at the origin $O$ with side of length $L(t)$ which solves the following degenerate ordinary differential equation

\begin{equation}
\dot{L}(t)=-{1\over \gamma}\left\lfloor\frac{\alpha_{hex}\gamma}{L(t)}\right\rfloor,
\label{unitareg2}
\end{equation}
until the extinction time.

We may compare this equation with the evolution law of the same initial sets by crystalline curvature with natural mobility that, as already remarked in the Introduction, can be defined independently of the Almgren, Taylor and Wang approach. The equation for $s(t)$ the apothem of $E(t)$ is given by (see \eqref{equamotions})

\begin{equation*}
\dot{s}(t)=-\varphi_{hex}(n_i)\kappa(t),
\end{equation*}
where $\varphi_{hex}(n_i)=\frac{2}{\sqrt{3}}$ for all $n_i\in\mathcal{N}$ and $\kappa(t)=\Lambda(n_i)/L(t)=\frac{4}{3L(t)}$ is the crystalline curvature. As a consequence, the corresponding law for the side length is 

\begin{equation}
\dot{L}(t)=-\alpha_{hex}\frac{1}{L(t)},
\label{mobreg}
\end{equation}
\\
showing that the limit evolution (\ref{unitareg2}), which is nontrivial only for sufficiently small Wulff shapes $E_0$, is slower than 
the corresponding crystalline evolutions. Moreover,

\begin{equation*}
\lim_{\gamma\to+\infty}-{1\over \gamma}\left\lfloor\frac{\alpha_{hex}\gamma}{L(t)}\right\rfloor=-\alpha_{hex}\frac{1}{L(t)}.
\end{equation*}
\end{oss}

\bigskip

\noindent
{\bf Acknowledgments.} I would like to thank Marco Cicalese for suggesting this problem and for his advices. 
I gratefully acknowledge the hospitality of the Technische Universit\"at  M\"unchen, where a part of this work has been carried out, and the financial support of the TU Foundation Fellowship.

\end{document}